\documentclass[journal]{IEEEtran}
\usepackage{CJK}
\usepackage{amssymb}
\usepackage{epstopdf}
\usepackage{graphicx}

\newtheorem{thm}{Theorem}
\newtheorem{lem}{Lemma}
\newtheorem{cor}{Corollary}

\newtheorem{rmk}{Remark}
\newtheorem{defi}{Definition}
\newtheorem{ass}{Assumption}

\newcommand {\emptycomment}[1]{}

\newcommand{\be }{\begin{equation}}
\newcommand{\ee }{\end{equation}}

\newcommand{\pf}{\emph{Proof: }}

\newcommand{\ttt}{\theta}
\newcommand{\ww}{\widetilde}

\newcommand{\f}{\frac}

\newcommand{\aaa}{\alpha}
\newcommand{\eqref}{\ref}

\newcommand{\nn}{\langle}
\newcommand{\mm}{\rangle}
\newcommand{\nono}{\nonumber}

\begin{document}
\begin{CJK*}{GBK}{song}

\title{Distributed Randomized Gradient-Free Mirror Descent Algorithm for Constrained Optimization}
\author{Zhan Yu, Daniel W. C. Ho,~\IEEEmembership{Fellow, IEEE}, Deming Yuan,~\IEEEmembership{Member, IEEE}
\thanks{Received November 6, 2018.
This work was partially supported by Research grants Council of the Hong Kong Special Administrative Region, China (CityU 11200717) and
CityU Strategic Research Grant 7005029.

Z. Yu and D. W. C. Ho are with the Department of Mathematics, City University of Hong Kong, Hong Kong  (e-mail: mathyuzhan@gmail.com;
madaniel@cityu.edu.hk).

D. Yuan is with the School of Automation, Nanjing University of Science and Technology, Nanjing 210094, China (e-mail: dmyuan1012@gmail.com)}
}
\maketitle

\begin{abstract}
This paper is concerned with multi-agent optimization problem. A distributed randomized gradient-free mirror descent (DRGFMD) method is developed by introducing a randomized gradient-free oracle in the mirror descent scheme where the non-Euclidean Bregman divergence is used. The classical gradient descent method is generalized without using subgradient information of objective functions. The proposed algorithm is the first distributed non-Euclidean zeroth-order method which achieves an $O(\f{1}{\sqrt{T}})$ convergence rate, recovering the best known optimal rate of distributed compact constrained convex optimization. Also, the DRGFMD algorithm achieves an $O(\f{\ln T}{T})$ convergence rate for the strongly convex constrained optimization case. The rate matches the best known non-compact constraint result. Moreover, a decentralized reciprocal weighted average approximating sequence is investigated and first used in distributed algorithm. A class of convergence rates are also achieved for the algorithm with weighted averaging (DRGFMD-WA). The technique on constructing the decentralized weighted average sequence provides new insight in searching for minimizers in distributed algorithms.
\end{abstract}



\IEEEpeerreviewmaketitle

\section{Introduction}

In recent years, distributed convex optimization over multi-agent network has played an important role in both theoretical and practical aspects. Early work mainly focuses on the research of minimizing a smooth function known to several agents (see, \cite{b1}-\cite{b3}). During these ten years, research has turned to the problem of minimizing a sum of locally convex objective functions distributed over a time-varying directed network via cooperation of agents(see, \cite{pd1}, \cite{s1}, \cite{s2}, \cite{directed}, \cite{sto}-\cite{g3}). The problem appears in diverse areas of science and engineering frequently. The seminal distributed method to solve the problem is by adopting distributed subgradient approach (see, \cite{s1}). Then, several methods applying distributed stochastic techniques to convex optimization take shape gradually (see, e.g. \cite{s2}, \cite{direct1}, \cite{direct2}). The stochastic subgradient method mainly considers stochastic disturbance on subgradient. The boundedness assumption of the stochastic gradient or subgradient and the consensus property among agents are two essential objects to ensure convergence of the algorithm. Moreover, the consensus property always relies on the topology of the network and the structure of the algorithm. In \cite{ss1}, the dual averaging structure is utilized for distributed optimization with an elegant probability approach. The distributed primal-dual algorithm has also been investigated recently in different directions (see, e.g. \cite{pd1}, \cite{pd3}). Aside from the methods above, the mirror descent, developed by Nemirovski and Yudin (see, \cite{mdmd}), attracts more and more research interest during these ten years. Thus, the mirror descent has been shown to be an efficient tool in several areas like large scale machine learning and sensor network. In \cite{mdc}, continuous-time stochastic mirror descent is investigated by Raginsky and Bouvrie by using Ito theory. In what follows, several discrete-time centralized and decentralized stochastic mirror descent methods are established quickly. Nedic and Lee consider a centralized mirror descent method for both convex and strongly convex optimization in \cite{md1}. Recently, Yuan et al. further develops an epoch type mirror descent method for strongly convex compact constrained optimization in \cite{md3}. Beside the convex optimization, several mirror descent methods have been established in non-convex optimization as well.

The approach in this paper to solve distributed nonsmooth optimization problem is related to randomized gradient-free method and classical mirror descent algorithm. Seminal randomized gradient-free techniques have been studied by Nesterov by considering different level of smoothness of functions (see, \cite{lemma tidu}) in non-distributed setting. After that, several researches on gradient-free method appear in distributed setting. Wang et al. \cite{g1} investigates a gradient-free method by taking the stochastic noise among agents into account. Sahu et al. \cite{g2} uses Kiefer-Wolfowitz gradient-free technique to solve a smooth optimization problem. Wang et al. \cite{g3} further extends the Kiefer-Wolfowitz method to distributed case. Hajinezhad et al. \cite{g4} establishes a gradient-free method to solve a nonconvex optimization problem via a primal-dual based framework. However, the existing distributed gradient-free convex optimization methods are all Euclidean projection based. This fact motivates us to consider a question: Is it possible to develop a distributed non-Euclidean gradient-free algorithm and obtain its convergence rate? In this paper, an answer is given in the affirmative. In this paper, the randomized gradient-free method is further developed by extending it to Bregman non-Euclidean framework and investigating the convergence rate of the proposed distributed algorithm. In distributed optimization problem over time-varying network, it is often necessary to build up consensus among agents to explore convergent behavior of the algorithm. Therefore an estimate of the expected disagreement among agents is given first. Then, the main results and corresponding convergence rates are established. In this work, we mention two common assumptions on objective functions and constraint set as follow: (1) Convex objective functions and compact convex constraint set; (2) Strongly convex objective functions and closed convex constraint set.  Some recent work indicates that, by choosing suitable stepsizes, assumption (1) and (2) often connect with a convergence rate of $O(\f{(\ln T)^{\alpha}}{T^{\beta}})$ type ($T$ is the number of iterations, $\alpha$ is nonnegative real number and $\beta$ is positive real number) (see, e.g. \cite{strongc}, \cite{sto}). The paper will analyze the proposed algorithm comprehensively by considering both assumptions. Moreover, based on the idea in \cite{md1} for centralized case, a decentralized weighted average approximating sequence is implemented in DRGFMD algorithm and several corresponding convergence rates are achieved.

The main theoretical contributions of this paper are summarized as follows:

\textbf{(i)} A decentralized zeroth-order (gradient-free) algorithm is proposed. The significance is that the algorithm carries Bregman non-Euclidean structure for solving the distributed convex optimization over time-varying network. As a result, the classical distributed zeroth-order projection algorithm is generalized to the non-Euclidean circumstance. Accordingly, the iteration procedure of the DRGFMD algorithm can provide a better reflection of the geometry of the convex optimization problem. The proposed algorithm is also operated in constrained domain with no smoothness requirement on objective functions. Meanwhile, different from existing mirror descent methods, the proposed algorithm relies on the random gradient-free oracles instead of knowing direct information on subgradients. Therefore, the algorithm becomes more flexible and efficiently applicable to areas like large scale machine learning and wireless sensor network where the subgradients of corresponding objective functions are infeasible or costly to evaluate.

\textbf{(ii)} To the best of our knowledge, in contrast to the existing methods, the proposed method is the first distributed zeroth-order (gradient-free) non-Euclidean method applied to convex and strongly convex optimization problems. In addition, a comprehensive analysis on DRGFMD is given under several conditions. For a convex optimization with compact constraint set, the proposed algorithms achieve the convergence rate of $O(\f{1}{\sqrt{T}})$. Thus, it recovers the best known convergence rate of the centralized compact constrained stochastic mirror descent algorithm in \cite{md1} and extends it to distributed situation.  In what follows,  an $O(\f{\ln T}{T})$ rate is obtained for strongly convex optimization problem with constraint set not assumed to be compact, extending the convergence results by Tsianos and Rabbat \cite{strongc} to non-Euclidean distributed situation. This is the first distributed zeroth-order non-Euclidean method to achieve it.

\textbf{(iii)} The paper investigates the reciprocal decentralized weighted average approximating sequence via the DRGFMD algorithm. Hence, it achieves a class of convergence rates for convex and strongly convex optimization problems. This is also the first distributed method that applies the reciprocal weighted average approximating sequence, therefore, the paper gives a future research direction on different types of decentralized weighted average sequences. It also provides a possibility to improve convergence rate in other distributed optimization algorithms.

\textbf{Notation:} Denote the n-dimension Euclidean space by $\mathbb R^n$, and the set of positive real numbers by $\mathbb R^{+}$. For a vector $v\in \mathbb R^n$, use $\|v\|$ to denote its Euclidean norm and $[v]_k$ to denote its $k$th component.  The inner  product of two vectors $a$, $b$ is denoted by $\nn a, b\mm$. For a matrix $M\in\mathbb R^{n\times n}$, denote the element in $i$th row and $j$th column by $[M]_{ij}$, denote the transpose of $M$ by $M^T$.  A function $f$ is $L$-Lipschitz on convex domain $X$ with respect to $\|\cdot\|$ if $|f(x)-f(y)|\leq L\|x-y\|$ holds for any $x, y\in X$.  A function $f$ is $\sigma_{f}$-strongly convex over domain $X$ if for any $x,y\in X$ and $\theta\in [0,1]$, $f(\ttt x+(1-\ttt)y)\leq\ttt f(x)+(1-\ttt)f(y)-\f{\sigma_{f}\ttt(1-\ttt)}{2}\|x-y\|^2$. Denote the gradient operator by $\nabla$, when $f$ is differentiable, the $\sigma_f$-strongly convex inequality above is equivalent to $f(x)\geq f(y)+\nn\nabla f(y),x-y\mm+\f{\sigma_{f}}{2}\|x-y\|^2$.
For two functions $f$ and $g$,  write $f(n)= O(g(n))$ if there exist $N<\infty$ and positive constant $C<\infty$ such that $f(n)\leq Cg(n)$ for $n\geq N$. For a random variable $X$, use $\mathbb E[X]$ to denote its expected value.

\section{Problem Setting}

In this paper, a time-varying multi-agent network is considered and the agents are indexed by $i=1,2,..., N$. The communication topology among agents is modeled as a directed graph $G_t=(V,E^t, P^t)$ in which $V=\{1,2,..., N\}$ is the node set, $E^t$ is the set of edges at time $t$, and $P^t$ is the communication matrix at time $t$. The goal of this work is to establish the distributed algorithms and convergence rate for the following distributed convex constrained optimization problem
\begin{eqnarray} \label{problem}
 \min_{x\in X} f(x)=\f{1}{N}\sum_{i=1}^Nf_i(x).
\end{eqnarray}
In (\eqref{problem}), $x\in X\subset\mathbb R^n$ is a global decision vector, $X$ is a nonempty convex domain. $f_i: \mathbb R^n\to \mathbb R$ is the convex objective function known only at the $i$th agent. $f_i$ is $L_i$-Lipschitz continuous over $X$.  In \textbf{Sections IV A} and \textbf{V A}, it is assumed that that each $f_i: \mathbb R^n\to \mathbb R$ are convex, not necessarily strongly convex, and $X$ is compact convex. In \textbf{Sections IV B} and \textbf{V B}, each $f_i$ is assumed to be strongly convex, and $X$ is closed convex, not necessarily compact. Denote the optimal point of the optimization problem by $x^{\star}$. In this paper, all the objective functions are supposed to be nonsmooth. Meanwhile, the situation when all the gradients or subgradients of the objective function can not be evaluated efficiently often appears. Thus, in this paper, the gradient-free technique is utilized to overcome this difficulty.  The smoothing function for objective function is introduced as a convolution of objective function $f_i$ and Gaussian kernel as follow,
\begin{eqnarray}
 f_{\mu_i}(x)=\f{1}{\kappa}\int_{\mathbb{R}^n}f_i(x+\mu_i\xi)e^{-\f{1}{2}\|\xi\|^2}d\xi,  \label{convolution}
\end{eqnarray}
in which $\kappa=(2\pi)^{\f{n}{2}}$ and $\mu_i\geq0$ is the smoothing parameter. Accordingly, the smoothing function of $f$ is denoted by
\begin{eqnarray}
  f_\mu(x)=\f{1}{N}\sum_{i=1}^Nf_{\mu_i}(x).  \label{defi}
\end{eqnarray}

In our non-Euclidean optimization algorithm, the Bregman divergence is considered as a non-Euclidean distance instead of the classical Euclidean distance employed by classical distributed gradient descent algorithms. The definition of the Bregman divergence is given below.
\begin{defi}
Let $\phi$ be a strongly convex differentiable function. The Bregman divergence between $x$ and $y$ induced by $\phi$ is denoted by $D_{\phi}(x,y)$ and given by $D_{\phi}(x,y)=\phi(x)-\phi(y)-\nn\nabla\phi(y), x-y\mm$.
\end{defi}
A basic result of  Bregman divergence is listed in the following lemma, the result will be used in subsequent analysis.  The proof follows from the definition directly.
\begin{lem}
The Bregman divergence satisfies the three-point identity $\nn\nabla\phi(x)-\nabla\phi(y),y-z\mm=D_{\phi}(z,x)-D_{\phi}(z,y)-D_{\phi}(y,x)$
for all $x, y, z\in X$.
\end{lem}

Some connections between the distance generating function $\phi$ and the properties of Bregman divergence will be described. The following assumption is made.

\begin{ass}\label{as1}
The distance generating function $\phi$ of Bregman divergence is three times continuously differentiable and $\sigma_{\phi}$-strongly convex with $\sigma>0$. The Hessian matrix $H=\nabla^2\phi$ and $H(y)+\nabla H(y)(y-x)$ are all positive semidefinite for any $x,y\in X$.
\end{ass}

Under Assumption \ref{as1} above, a direct consequence is the relation between Bregman divergence and the classical Euclidean distance: $D_{\phi}(x,y)\geq \f{\sigma_{\phi}}{2}\|x-y\|^2$. Another important consequence obtained from Assumption \ref{as1} is the separate convexity of Bregman divergence $D_{\phi}(x,y)$: $D_{\phi}(x, \sum_{j=1}^N\theta_jy_j)\leq\sum_{j=1}^N\theta_jD_{\phi}(x,y_j)$
in which $\sum_{j=1}^N\theta_j=1$, $\theta_j\geq 0$.

In what follows, the standard assumption on the graph $G_t=(V,E^t, P^t)$ is made.

\begin{ass}\label{as2} Communication matrix $P^t$ is a doubly stochastic matrix, $i.e.$, $\sum_{i=1}^N[P^t]_{ij}=1$ and $\sum_{j=1}^N[P^t]_{ij}=1$ for any $i$ and $j$. There exists some positive integer $B$ such that the graph $(V, \bigcup_{s=0}^{B-1}E^{t+s})$ is strongly connected for any $t$. There exists a scalar $0<\zeta<1$ such that $[P^t]_{ii}\geq \zeta$ for all $i$ and $t$, and $[P^t]_{ij}\geq\zeta$ if $(j,i)\in E^t$.
\end{ass}

Denote the transition matrices by $P(t,s)=P^tP^{t-1}\cdots P^s, t\geq s\geq0$, an important consequence about the transition matrices is listed in the following lemma. The result will be essential in subsequent analysis.
\begin{lem}\label{1n}
\cite{s1}Let Assumption \ref{as2} hold, then for all $i,j\in V$ and all $t,s$ satisfying $t\geq s\geq 0$
\begin{eqnarray}
|[P(t,s)]_{ij}-\f{1}{N}|\leq \Gamma\gamma^{t-s},
\end{eqnarray}
in which $\Gamma=(1-\f{\zeta}{4N^2})^{-2}$ and $\gamma=(1-\f{\zeta}{4N^2})^{\f{1}{B}}$.
\end{lem}
The following auxiliary lemma for sequence is basic for later use.
\begin{lem}\label{sequence}
\cite{s2}Given a positive sequence $\{\beta_k\}_{k\geq0}$ with $\lim_{k\to \infty}\beta_k=\ww{\beta}$ and $\gamma\in(0,1)$, the following holds,
\begin{eqnarray}
\lim_{t\to\infty}\sum_{k=1}^{t-1}\gamma^{t-k}\beta_{k-1}=\f{\gamma}{1-\gamma}\ww{\beta}.
\end{eqnarray}
\end{lem}

\section{The algorithm and preliminary for convergence theorem}
Let $x_i^t$ be the estimate of agent $i$ at step $t$. The distributed randomized gradient-free mirror descent (DRGFMD) algorithm is designed as
\begin{eqnarray}
 y_i^t&=&\sum_{j=1}^N[P^t]_{ij}x_j^t, \label{ag1} \\
  x_i^{t+1}&=&\arg\min_{x\in X}\{\aaa_t\nn g_{\mu_i}(z_i^t),x\mm+D_{\phi}(x,y_i^t)\}, \label{ag2}
\end{eqnarray}
in which $z_i^t$ denotes $y_i^t$ or $x_i^t$ for any $i\in V$ and $t\in\mathbb N$. Here,
\begin{eqnarray}
g_{\mu_i}(z_i^t)=\f{f_i(z_i^t+\mu_i\xi_i^t)-f_i(z_i^t)}{\mu_i}\xi_i^t
\end{eqnarray}
is the random gradient-free oracle. $\xi_{i}^t$ is a random vector sequence that is locally generated in an i.i.d distributed manner according to the Gaussian distribution for each agent $i\in V$. The information communicating behavior among agent $i$ and its neighbors at $t$ step is described in (\eqref{ag1}). In this step, agent $i$ receives estimates $x_j^t$ from its neighbors $j\in N_i^t$, computes a weighted average on all the received estimates to get a new state variable $y_i^t$. Then in (\eqref{ag2}), the algorithm updates locally via a gradient-free approach which is built on the mirror descent scheme. In following proofs in \textbf{Section III}, \textbf{Section IV} and \textbf{Section V}, without loss of generality, it is assumed that $x_i^0=0$. The algorithm starts at $x_i^0$, the estimates of agents used to construct the approximating sequence start from $x_i^1$.
\begin{rmk}
In this algorithm, the classical gradient descent method is generalized in two aspects. On one hand, Bregman divergence is utilized instead of the classical Euclidean norm, leading to the non-Euclidean projection feature of the proposed algorithms. The non-Euclidean form of the algorithm can better reflect the implicit geometry characteristic of the distributed optimization problem. Several classical choices of distance generating function $\phi$ make the algorithm more efficient than other algorithms when applied to many areas like large-scale machine learning and wireless sensor networks. On the other hand, the random gradient-free oracle in the algorithm enables us to apply the algorithm to situation where the subgradient of objective function is hard to achieve, thus overcomes the tough environment on subgradient successfully.
\end{rmk}

Now some preliminaries are made for proving the main results. Firstly several important estimates on the random oracle are made in following lemma. Define the $\sigma$-field generated by the history of the random variables to step $t-1$ in the form: $F_t=\{(x_i^0, i=1, 2,..., N); (\xi_i^s, i=1, 2,..., N); 1\leq s\leq t-1\}$ with $F_0=\{x_i^0,i=1, 2,..., N\}$. Then the following lemma holds.
\begin{lem}\label{lip}\cite{gf}
Let $\widehat{L}=\max_{i}L_i$, for $x\in X$ the following holds.

(a)$f(x)$ satisfies $f_i(x)\leq f_{\mu_i}(x)\leq f_i(x)+\sqrt{n}\mu_i\widehat{L}$.

(b)$f_{\mu_i}(x)$ is convex, differentiable and the following relation with the oracle $g_{\mu_i}(z_i^t)$ holds: $\mathbb E[g_{\mu_i}(z_i^t)|F_{t}]=\nabla f_{\mu_i}(z_i^t)$.

(c)The random gradient-free oracle $g_{\mu_i}(z_i^t)$ satisfies $\mathbb E[\|g_{\mu_i}(z_i^t)\|^2|F_{t}]\leq(n+4)^2\widehat{L}^2$.
\end{lem}
By using Cauchy inequality and Minkowski inequality, two direct results can be gotten from Lemma \ref{lip}:
\begin{eqnarray}
\nonumber \mathbb E[\|g_{\mu_i}(z_i^t)\||F_{t}]\leq(n+4)\widehat{L}, \forall i\in V,
\end{eqnarray}
and
\begin{eqnarray}
\|\nabla f_{\mu_i}(z_i^t)\|\leq (n+4)\widehat{L}, \forall i\in V. \label{kexi}
\end{eqnarray}
Denote the Bregman projection error by
\begin{eqnarray}
e_i^t=x_i^{t+1}-y_i^t,\label{error}
\end{eqnarray}
the following lemma gives an upper bound estimate of the error. The estimate is necessary to obtain corresponding consensus property of the estimates for each agent, which guarantees the convergence of the algorithm.
\begin{lem}\label{errorbd}
Under Assumption \ref{as1}, for the DRGFMD algorithm, let the Bregman projection error for agent $i$ be defined as (\eqref{error}). Then for any $i\in V$ and $t\geq 0$, the following holds,
\begin{eqnarray}
\nonumber \mathbb E[\|e_i^t\|]\leq \f{(n+4)\widehat{L}}{\sigma_{\phi}}\aaa_t.
\end{eqnarray}
\end{lem}

\pf The first-order optimality of $x_i^{t+1}$ implies
\begin{eqnarray}\label{opt}
\nn \aaa_t g_{\mu_i}(z_i^t)+\nabla\phi(x_i^{t+1})-\nabla\phi(y_i^t),x-x_i^{t+1}\mm\geq0, \ \forall x\in X.\label{optcon}
\end{eqnarray}
By setting $x=y_i^t$ in (\eqref{optcon}) and using the $\sigma_{\phi}$-strongly convexity of $\phi$, it follows that
\begin{eqnarray}
 \nonumber \nn\aaa_t g_{\mu_i}(z_i^t),y_i^t-x_i^{t+1}\mm&\geq&\nn\nabla\phi(y_i^t)-\nabla\phi(x_i^{t+1}),y_i^t-x_i^{t+1}\mm\\
\nonumber  &\geq&\sigma_{\phi}\|y_i^t-x_i^{t+1}\|^2.
\end{eqnarray}
By using Cauchy inequality to the left hand side of the above inequality, the following holds,
\begin{eqnarray}
\nonumber \aaa_t\|g_{\mu_i}(z_i^t)\|\cdot\|y_i^t-x_i^{t+1}\|\geq\ \sigma_{\phi}\|y_i^t-x_i^{t+1}\|^2.
\end{eqnarray}
Dividing by $\|y_i^t-x_i^{t+1}\|$ on both sides yields
\begin{eqnarray}
\nonumber \|y_i^t-x_i^{t+1}\|\leq\f{1}{\sigma_{\phi}}\|g_{\mu_i}(z_i^t)\|\cdot\aaa_t, \ \forall i\in V, \ t\geq 0 .
\end{eqnarray}
By taking the conditional expectation on $F_t$  and using Lemma \ref{lip}, the final result is obtained after taking the total expectation.

Next, the consensus result among agents is ready to be established. The average of all agents at step $t$ is introduced as follow:
\begin{eqnarray}
\overline{x}^{t}=\f{1}{N}\sum_{i=1}^{N}x_i^{t}. \label{aver}
\end{eqnarray}
\begin{lem}\label{consensus}
Under Assumption \ref{as2}, let $\{x_i^t\}_{t\geq 0}$ be the sequence generated by DRGFMD algorithm. Then for any non-increasing positive stepsizes $\aaa_t$ and any agent $i, j\in V$, the following estimate holds:
\begin{small}
\begin{eqnarray}
\nonumber \mathbb E[\|x_i^{t}-\overline{x}^{t}\|]\leq\f{N\Gamma(n+4)\widehat{L} }{\sigma_{\phi}}\sum_{k=1}^{t-1}\gamma^{t-k}\aaa_{k-1}+\f{2(n+4)\widehat{L}}{\sigma_{\phi}}\aaa_{t-1};\\
\nonumber \mathbb E[\|x_i^t-x_j^{t}\|]\leq\f{2N\Gamma(n+4)\widehat{L} }{\sigma_{\phi}}\sum_{k=1}^{t-1}\gamma^{t-k}\aaa_{k-1}+\f{4(n+4)\widehat{L}}{\sigma_{\phi}}\aaa_{t-1};\\
\nonumber \sum_{t=1}^T\sum_{i=1}^N\mathbb E[\|x_i^t-x_j^{t}\|]\leq (\f{2N^2\Gamma(n+4)\widehat{L}\gamma }{\sigma_{\phi}(1-\gamma)}+\f{4N(n+4)\widehat{L}}{\sigma_{\phi}})\sum_{t=0}^T\aaa_t.
\end{eqnarray}
\end{small}
\end{lem}
\pf By iterating recursively, the update $x_i^{t}$ can be expanded in the form
\begin{eqnarray}
\nonumber x_i^{t} &=& y_i^{t-1}+e_i^{t-1}\\
  &=& \sum_{j=1}^N[P(t-1,0)]_{ij}x_j^0+\sum_{k=1}^{t-1}\sum_{j=1}^N[P(t-1,k)]_{ij}e_j^{k-1}+e_i^{t-1}. \label{rec1}
\end{eqnarray}
Taking average on both sides and noting that $P^t$ is doubly stochastic, $\overline{x}^{t}$ can be written in the following form,
\begin{eqnarray}
 \overline{x}^{t}=\f{1}{N}\sum_{j=1}^Nx_j^0+\f{1}{N}\sum_{k=1}^t\sum_{j=1}^Ne_j^{k-1}. \label{rec2}
\end{eqnarray}
Combining (\eqref{rec1}) and (\eqref{rec2}) yields
\begin{eqnarray}
\nonumber &&\|x_i^{t}-\overline{x}^{t}\|\\
\nonumber&=&\|\sum_{j=1}^N([P(t-1,0)]_{ij}-\f{1}{N})x_j^0+\sum_{k=1}^{t-1}\sum_{j=1}^N([P(t-1,k)]_{ij}-\f{1}{N}) \\
\nonumber&&\cdot e_j^{k-1}+(e_i^{t-1}-\f{1}{N}\sum_{j=1}^Ne_j^{t-1})\|\\
 \nonumber &\leq& \sum_{j=1}^N|[P(t-1,0)]_{ij}-\f{1}{N}|\cdot\|x_j^0\|+\sum_{k=1}^{t-1}\sum_{j=1}^N|[P(t-1,k)]_{ij}\\
 &&-\f{1}{N}|\cdot \|e_j^{k-1}\|+\f{1}{N}\sum_{j=1}^N\|e_j^{t-1}\|+\|e_i^{t-1}\|.
\end{eqnarray}
Take total expectation on both sides of the inequality above and note that $x_j^0=0, j\in V$, the first consensus result follows from Lemma \ref{1n} and Lemma \ref{errorbd}. The second one is obtained directly by using the triangle inequality to the first one. Sum $\mathbb E[\|x_i^t-x_j^t\|]$ up over the indices from $t=1$ to $T$ and $i=1$ to $N$, it follows that
\begin{eqnarray}
\nonumber&&\sum_{t=1}^T\sum_{j=1}^N\mathbb E[\|x_i^t-x_j^t\|]\\
\nonumber&\leq&\f{2N^2\Gamma(n+4)\widehat{L} }{\sigma_{\phi}}\sum_{t=1}^T\sum_{k=1}^{t-1}\gamma^{t-k-1}\aaa_{k-1}+\f{4N(n+4)\widehat{L}}{\sigma_{\phi}}\sum_{t=1}^T\aaa_{t-1} \\
\nonumber&\leq&\f{2N^2\Gamma(n+4)\widehat{L} }{\sigma_{\phi}}\sum_{s=0}^{T-2}\gamma^s\sum_{t=0}^{T-2}\aaa_{t}+\f{4N(n+4)\widehat{L}}{\sigma_{\phi}}\sum_{t=0}^{T-1}\aaa_{t}  \\
\nonumber&\leq&\f{2N^2\Gamma(n+4)\widehat{L}}{\sigma_{\phi}}\cdot\f{1}{1-\gamma}\sum_{t=0}^T\aaa_t+\f{4N(n+4)\widehat{L}}{\sigma_{\phi}}\sum_{t=0}^T\aaa_t,
\end{eqnarray}
combining the two terms completes the proof.

The following result provides an essential estimate for the main convergence results, it will also be useful in following sections. This part is concluded by analyzing this inner product estimate result.
\begin{thm}\label{neiji}
Under Assumptions \ref{as1} and \ref{as2}, let $\{x_i^t\}_{t\geq 0}$ and $\{y_i^t\}_{t\geq 0}$ be the sequences generated by the DRGFMD algorithm, $\aaa_t$ is any non-increasing sequence, then for general convex domain $X$ and any $T\geq 1$, we have
\begin{eqnarray}
\nonumber \nn\aaa_t g_{\mu_i}(z_i^t), y_i^t-x^{\star}\mm&\leq& D_{\phi}(x^{\star},y_i^t)-D_{\phi}(x^{\star},x_i^{t+1}) \\
 &&+\f{\aaa_t^2}{2\sigma_{\phi}}\|g_{\mu_i}(z_i^t)\|^2,  \label{sides}
\end{eqnarray}
and
\begin{small}
\begin{eqnarray}
\nonumber &&\f{1}{NT}\sum_{t=1}^T\sum_{i=1}^N\mathbb E[\nn \nabla f_{\mu_i}(z_i^t), y_i^t-x^{\star}\mm] \\
\nonumber &\leq&\f{1}{NT}\sum_{i=1}^N\big[\f{1}{\aaa_1}\mathbb E [D_{\phi}(x^{\star},x_i^1)]   \\
\nonumber&&+\sum_{t=2}^T\mathbb E[D_{\phi}(x^{\star}, x_i^t)](\f{1}{\aaa_t}-\f{1}{\aaa_{t-1}})\\
&& -\f{1}{\aaa_T}\mathbb E[D_{\phi}(x^{\star}, x_i^{T+1})]\big]+\f{(n+4)^2\widehat{L}^2}{2\sigma_{\phi} T}\sum_{t=0}^T\aaa_t.   \label{inner1}
\end{eqnarray}
\end{small}
Moreover, when $X$ is a compact convex constraint set with a constant $d_{\phi}$ such that $\sup_{x,y\in X} D_{\phi}(x, y)\leq d_{\phi}^2$. Then the following holds,
\begin{eqnarray}
 \nono&&\f{1}{NT}\sum_{t=1}^T\sum_{i=1}^N\nn \nabla f_{\mu_i}(z_i^t), y_i^t-x^{\star}\mm\\
 &&\leq \f{d_{\phi}^2}{T\aaa_T}+\f{(n+4)^2\widehat{L}^2}{2\sigma_{\phi} T}\sum_{t=0}^T\aaa_t.  \label{inner2}
\end{eqnarray}

\end{thm}
\pf Set $x=x^{\star}$ in (\eqref{opt}) and use the Bregman three-point inequality, it can be obtained that
\begin{eqnarray}
 \nonumber&&\nn\aaa_t g_{\mu_i}(z_i^t), x_i^{t+1}-x^{\star}\mm \\
\nonumber &\leq&\nn\nabla\phi(y_i^t)-\nabla\phi(x_i^{t+1}), x_i^{t+1}-x^{\star}\mm\\
 \nonumber &=& D_{\phi}(x^{\star},y_i^t)-D_{\phi}(x^{\star},x_i^{t+1})-D_{\phi}(x_i^{t+1},y_i^t) \\
  &\leq&D_{\phi}(x^{\star},y_i^t)-D_{\phi}(x^{\star},x_i^{t+1})-\f{\sigma_{\phi}}{2}\|x_i^{t+1}-y_i^t\|^2,\label{nei1}
\end{eqnarray}
in which the second inequality follows from the definition of Bregman divergence and the $\sigma_{\phi}$-strongly convexity of $\phi$. On the other hand,
\begin{eqnarray}
\nonumber  &&\nn\aaa_t g_{\mu_i}(z_i^t), x_i^{t+1}-x^{\star}\mm\\
\nonumber&=&\nn\aaa_t g_{\mu_i}(z_i^t), x_i^{t+1}-y_i^t\mm+\nn\aaa_t g_{\mu_i}(z_i^t), y_i^t-x^{\star}\mm\\
  \nonumber&\geq&-\f{\aaa_t^2}{2\sigma_{\phi}}\|g_{\mu_i}(z_i^t)\|^2-\f{\sigma_{\phi}}{2}\|x_i^{t+1}-y_i^t\|^2\\
  &&+\nn\aaa_t g_{\mu_i}(z_i^t), y_i^t-x^{\star}\mm,\label{nei2}
\end{eqnarray}
where the second inequality follows from Fenchel inequality. Thus, (\eqref{nei1}) and (\eqref{nei2}) together imply (\eqref{sides}). Divide by $\aaa_t$ on both sides of (\eqref{sides}), it follows
\begin{eqnarray}
 \nono\nn g_{\mu_i}(z_i^t), y_i^t-x^{\star}\mm&\leq& \f{1}{\aaa_t}[D_{\phi}(x^{\star},y_i^t)-D_{\phi}(x^{\star},x_i^{t+1})]\\
 &&+\f{\aaa_t}{2\sigma_{\phi}}\|g_{\mu_i}(z_i^t)\|^2. \label{nei3}
\end{eqnarray}
Sum (\eqref{nei3}) over the indices from $t=1$ to $T$ and $i=1$ to $N$, it follows
\begin{small}
\begin{eqnarray}
 \nonumber&&\sum_{t=1}^T\sum_{i=1}^N\nn g_{\mu_i}(z_i^t),  y_i^t-x^{\star}\mm  \\
 \nonumber &\leq&\sum_{t=1}^T\f{1}{\aaa_t}[\sum_{i=1}^N\sum_{j=1}^N[P^t]_{ij}D_{\phi}(x^{\star}, x_j^t)-\sum_{i=1}^ND_{\phi}(x^{\star}, x_i^{t+1})]\\
 \nono&&+\sum_{t=1}^T\sum_{i=1}^N\f{\aaa_t}{2\sigma_{\phi}}\|g_{\mu_i}(z_i^t)\|^2 \\
\nonumber  &=&\sum_{t=1}^T\f{1}{\aaa_t}[\sum_{j=1}^ND_{\phi}(x^{\star}, x_j^t)-\sum_{i=1}^ND_{\phi}(x^{\star}, x_i^{t+1})] \\
\nono&&+\sum_{t=1}^T\sum_{i=1}^N\f{\aaa_t}{2\sigma_{\phi}}\|g_{\mu_i}(z_i^t)\|^2 \\
 \nonumber &=& \sum_{i=1}^N[\f{1}{\aaa_1}D_{\phi}(x^{\star},x_i^1)+\sum_{t=2}^TD_{\phi}(x^{\star}, x_i^t)(\f{1}{\aaa_t}-\f{1}{\aaa_{t-1}})\\
  &&-\f{1}{\aaa_T}D_{\phi}(x^{\star}, x_i^{T+1})]+\sum_{t=1}^T\sum_{i=1}^N\f{\aaa_t}{2\sigma_{\phi}}\|g_{\mu_i}(z_i^t)\|^2, \label{nei4}
\end{eqnarray}
\end{small}
where the first inequality follows from the separate convexity of Bregman divergence, first equality follows from doubly stochastic property of matrix $P^t$, and the second equality is as a result of rearranging terms.  Taking conditional expectation on both sides of (\eqref{nei4}) over $F_t$, then taking total expectation, using Lemma \ref{lip} and noting that $\sum_{t=1}^T\aaa_t\leq\sum_{t=0}^T\aaa_t$ yields (\eqref{inner1}). When $X$ is compact, divide by $NT$ on both sides of (\eqref{nei4}), note that $-\f{1}{\aaa_T}\mathbb E[D_{\phi}(x^{\star}, x_i^{T+1})]\leq0$, it follows
\begin{small}
\begin{eqnarray}
 \nono&&\f{1}{NT}\sum_{t=1}^T\sum_{i=1}^N\nn g_{\mu_i}(z_i^t), y_i^t-x^{\star}\mm\\
 \nono&\leq& \f{1}{NT}\sum_{i=1}^Nd_{\phi}^2[\f{1}{\aaa_1}+\sum_{t=2}^T(\f{1}{\aaa_t}-\f{1}{\aaa_{t-1}})] \\
 \nono&&+\f{1}{2\sigma_{\phi} NT}\sum_{t=1}^T\sum_{i=1}^N\|g_{\mu_i}(z_i^t)\|^2\aaa_t,
\end{eqnarray}
\end{small}
namely,
\begin{eqnarray}
\nono && \f{1}{NT}\sum_{t=1}^T\sum_{i=1}^N\nn g_{\mu_i}(z_i^t), y_i^t-x^{\star}\mm\leq \f{d_{\phi}^2}{T\aaa_T}\\
 &&+\f{1}{2\sigma_{\phi} NT}\sum_{t=1}^T\sum_{i=1}^N\|g_{\mu_i}(z_i^t)\|^2\aaa_t.\label{nei5}
\end{eqnarray}
Take conditional expectation on both sides of (\eqref{nei5}) over $F_t$ and use Lemma \ref{lip} again,  (\eqref{inner2}) is obtained.

\section{The DRGFMD algorithm with classical approximating sequence}
In this section, the DRGFMD algorithm with the classical approximating sequence
\begin{eqnarray}
\nonumber \widehat{x}_{l}^{T}=\f{1}{T}\sum_{t=1}^Tx_l^t
\end{eqnarray}
 is applied to both convex and strongly convex optimization problem. The convergence theorems and convergence rates for them will be established respectively. Several advantages of the algorithm will be discussed in detail.

\subsection{DRGFMD algorithm for compact constrained convex optimization}
The distributed randomized gradient-free mirror descent algorithm under compact convex constrained condition is studied in this part.
Equipped with Lemma \ref{consensus} and Theorem \ref{neiji}, the main convergence results are ready to be presented for the following DRGFMD algorithm:
\begin{small}
\begin{eqnarray}
 \nonumber y_i^t&=&\sum_{j=1}^N[P^t]_{ij}x_j^t, \\
\nonumber  x_i^{t+1}&=&\arg\min_{x\in X}\{\aaa_t\nn \f{f_i(y_i^t+\mu_i\xi_i^t)-f_i(y_i^t)}{\mu_i}\xi_i^t,x\mm+D_{\phi}(x,y_i^t)\}.
\end{eqnarray}
\end{small}

\begin{thm}\label{main}
Under Assumptions \ref{as1} and \ref{as2}, let $\{x_i^t\}_{t\geq 0}$ and $\{y_i^t\}_{t\geq 0}$ be the sequences generated by the DRGFMD algorithm with a non-increasing positive stepsize sequence $\aaa_t$. Then, for any optimal point $x^{\star}$ and $l\in V$, the following convergence result of $\widehat{x}_{l}^{T}$ for compact constrained convex optimization holds:
\begin{eqnarray}
&&\mathbb E[f(\widehat{x}_l^T)]-f(x^{\star})\leq B_1+B_2+B_3. \label{thm2}
\end{eqnarray}
\begin{eqnarray}
\nonumber  B_1=\sqrt{n}\widehat{L}\cdot\f{1}{N}\sum_{i=1}^N\mu_i; \ \ \ \  B_2=\f{d_{\phi}^2}{T\aaa_T};
\end{eqnarray}
\begin{eqnarray}
\nonumber  B_3=[\f{9(n+4)^2\widehat{L}^2}{2\sigma_{\phi} }+\f{2N(n+4)^2\Gamma\widehat{L}^2}{\sigma_{\phi} (1-\gamma)}]\cdot\f{1}{T}\sum_{t=0}^T\aaa_t.
\end{eqnarray}
\end{thm}
\pf According to the convexity of $f_{\mu_i}(x)$ at point $y_i^t$ (convexity of $f_{\mu_i}$ follows from the convexity of $f_i$),
\begin{eqnarray}
 \f{1}{N}\sum_{i=1}^N\nn \nabla f_{\mu_i}(y_i^t), y_i^t-x^{\star}\mm\geq\f{1}{N}\sum_{i=1}^N( f_{\mu_i}(y_i^t)-f_{\mu_i}(x^{\star})). \label{neifmu}
\end{eqnarray}
For any $l\in V$, the following estimate of $\f{1}{N}\sum_{i=1}^Nf_{\mu_i}(y_i^t)$ in (\eqref{neifmu}) is made as follow,
\begin{eqnarray}
\nonumber \f{1}{N}\sum_{i=1}^Nf_{\mu_i}(y_i^t)&=&\f{1}{N}\sum_{i=1}^N[f_{\mu_i}(x_l^t)+f_{\mu_i}(y_i^t)-f_{\mu_i}(x_l^t)]\\
 \nonumber &\geq&f_{\mu}(x_l^t)-\f{1}{N}\sum_{i=1}^N\|\nabla f_{\mu_i}(y_i^t)\|\cdot\|y_i^t-x_l^t\| \\
  &\geq&f_{\mu}(x_l^t)-\f{(n+4)\widehat{L}}{N}\sum_{i=1}^N\|y_i^t-x_l^t\|,
\end{eqnarray}
in which the first inequality follows from the convexity of $f_{\mu_i}$ and Cauchy inequality, the second inequality follows from (\eqref{kexi}). Notice that $\|y_i^t-x_l^t\|=\|\sum_{j=1}^N[P^t]_{ij}x_j^t-x_l^t\|\leq\sum_{j=1}^N[P^t]_{ij}\|x_j^t-x_l^t\|$, it follows that
\begin{eqnarray}
 \f{1}{N}\sum_{i=1}^Nf_{\mu_i}(y_i^t)\geq f_{\mu}(x_l^t)-\f{(n+4)\widehat{L}}{N}\sum_{j=1}^N\|x_j^t-x_l^t\|, \label{fmu1}
\end{eqnarray}
in which the equality follows from the doubly stochastic property of matrix $P^t$. Substituting (\eqref{fmu1}) into (\eqref{neifmu}) yields
\begin{eqnarray}
\nonumber &&\f{1}{N}\sum_{i=1}^N\nn \nabla f_{\mu_i}(y_i^t), y_i^t-x^{\star}\mm \\
\nono&\geq& \f{1}{N}\sum_{i=1}^Nf_{\mu_i}(y_i^t)-\f{1}{N}\sum_{i=1}^Nf_{\mu_i}(x^{\star}) \\
\nonumber &\geq&f(x_l^t)-f_{\mu}(x^{\star})-\f{(n+4)\widehat{L}}{N}\sum_{j=1}^N\|x_j^t-x_l^t\| \\
\nono&\geq& f(x_l^t)-f(x^{\star})-\sqrt{n}\widehat{L}\cdot\f{1}{N}\sum_{i=1}^N\mu_i  \\
&&-\f{(n+4)\widehat{L}}{N}\sum_{j=1}^N\|x_j^t-x_l^t\|, \label{fmu2}
\end{eqnarray}
in which the third inequality follows from Lemma \ref{lip}. Summing up both sides of (\eqref{fmu2}) from $t=1$ to $T$ and dividing by $T$, then taking total expectation on both sides implies
\begin{eqnarray}
\nonumber &&\f{1}{NT}\sum_{t=1}^T\sum_{i=1}^N\mathbb E[\nn \nabla f_{\mu_i}(y_i^t), y_i^t-x^{\star}\mm] \\
\nono &&\geq \f{1}{T}\sum_{t=1}^T\mathbb E[f(x_l^t)]-f(x^{\star})-\sqrt{n}\widehat{L}\cdot\f{1}{N}\sum_{i=1}^N\mu_i \\
 &&-\f{(n+4)\widehat{L}}{NT}\sum_{t=1}^T\sum_{j=1}^N\mathbb E[\|x_j^t-x_l^t\|]. \label{fmu3}
\end{eqnarray}
Rearrange terms, combine (\eqref{fmu3}) and Theorem \ref{neiji}, use the convexity of $f$, it can be obtained that
\begin{eqnarray}
\nonumber&& \mathbb E[f(\widehat{x}_l^T)]-f(x^{\star})\\
\nono&&\leq \sqrt{n}\widehat{L}\cdot\f{1}{N}\sum_{i=1}^N\mu_i+\f{d_{\phi}^2}{T\aaa_T}+\f{(n+4)^2\widehat{L}^2}{2\sigma_{\phi} T}\sum_{t=0}^T\aaa_t\\
&&+\f{(n+4)\widehat{L}}{NT}\sum_{t=1}^T\sum_{j=1}^N\mathbb E[\|x_j^t-x_l^t\|]. \label{xiao}
\end{eqnarray}
The desired convergence result follows by applying Lemma \ref{consensus} to the last term of (\eqref{xiao}).

\begin{rmk}
Theorem \ref{main} indicates the convergence property of the local sequence $x_i^t$ via the average vector $\widehat{x}_l^T=\f{1}{T}\sum_{t=1}^Tx_l^{t}$ at each node $i\in V$. In fact, by taking sequence $\widehat{\overline{x}^T}=\f{1}{T}\sum_{t=1}^T\overline{x}^{t}$, the algorithm can also generate an approximating convergence sequence $\widehat{\overline{x}^T}$ with constant $\f{9(n+4)^2\widehat{L}^2}{2\sigma_{\phi} }+\f{2N(n+4)^2\Gamma\widehat{L}^2}{\sigma_{\phi} (1-\gamma)}$ replaced by $\f{5(n+4)^2\widehat{L}^2}{2\sigma_{\phi} }+\f{N(n+4)^2\Gamma\widehat{L}^2}{\sigma_{\phi} (1-\gamma)}$. In order to highlight the local characteristic of convergence sequence for each node $i\in V$,  $\widehat{x}_l^T$ is used instead of $\widehat{\overline{x}^T}$ as the approximating sequence.
\end{rmk}

\begin{rmk}
Theorem \ref{main} shows that the expected convergence error in (\eqref{thm2}) of the proposed algorithm is upper bounded by three parts. The first part $B_1$ is the smoothing parameters as a penalty of using the gradient-free oracle instead of the true gradient information. By selecting arbitrary small parameter $\mu_i$, the optimality gap and convergence error are reduced. The second part $B_2$ describes the influence of the structure of the domain and  non-Euclidean structure of the underlying metric. The third part $B_3$ mainly describes the influence of the topology of the network under the non-Euclidean structure. We can select appropriate diminishing stepsize $\aaa_t$ to handle convergence rate.
\end{rmk}


\begin{thm}\label{crate1}
Under Assumptions \ref{as1} and \ref{as2}, let $\{x_i^t\}_{t\geq 0}$ and $\{y_i^t\}_{t\geq 0}$ be the sequences generated by the DRGFMD algorithm. Select $\aaa_t=\f{\eta}{\sqrt{t+1}}$ ($\eta$ is a positive parameter). Then for any $l\in V$, any optimal point $x^{\star}$ and $T\geq1$, $\widehat{x}_l^T$ converges asymptotically to the approximate optimal solution in following sense:
\begin{eqnarray}
\mathbb E[f(\widehat{x}_l^T)]-f(x^{\star})\leq \sqrt{n}\widehat{L}\cdot\f{1}{N}\sum_{i=1}^N\mu_i+\f{C_1}{\sqrt{T}}, \label{rate}
\end{eqnarray}
$C_1$ is a constant as follow,
\begin{eqnarray}
\nono C_1=\f{\sqrt{2}d_{\phi}^2}{\eta}+\f{9\sqrt{2}(n+4)^2\widehat{L}\eta}{\sigma_{\phi} }+\f{4\sqrt{2}N(n+4)^2\Gamma\widehat{L}^2\eta\gamma}{\sigma_{\phi} (1-\gamma)}.
\end{eqnarray}
\end{thm}
\pf Since the  stepsize $\aaa_t=\f{\eta}{\sqrt{t+1}}$, it follows
\begin{eqnarray}
\nono&&\sum_{t=0}^T\aaa_t=\eta+\sum_{t=1}^T\f{\eta}{\sqrt{t+1}}\leq\eta+\int_{0}^T\f{\eta}{\sqrt{t+1}}dt \\
&&\leq 2\eta\sqrt{T+1}-\eta\leq2\eta\sqrt{T+1}.
\end{eqnarray}
Substituting it into Theorem \ref{main} and noticing that $\sqrt{T+1}\leq\sqrt{2T}$ for $T\geq1$ yields (\eqref{rate}).

\begin{rmk}
Theorem \ref{crate1} indicates that the DRGFMD algorithm converges at an $O(\f{1}{\sqrt{T}})$ rate in distributed compact constrained convex optimization problem. The rate matches the best known compact constraint convergence rate of distributed subgradient algorithm in Yuan et al. (\cite{sto}). Moreover, the paper extends the methods to a non-Euclidean distributed scenario. To the best of our knowledge, the proposed algorithm is the first distributed algorithm which makes use of gradient-free technique to implement a distributed mirror descent algorithm. Moreover, the convergence rate is obtained without any smoothness assumptions on objective functions, making the algorithm applicable to more extensive areas in science and engineering.
\end{rmk}

\begin{rmk}
In (\eqref{rate}), it provides intrinsic information of the influence of the dimension of agents' estimates $n$ on convergence error of DRGFMD, which shows that the convergence gets slower when $n$ becomes larger. This phenomenon will be illustrated in the simulations.
\end{rmk}

\begin{rmk}
In fact, if the estimator $g_{\mu_i}(x_i^t)$ is utilized in the DRGFMD algorithm, an $O(\f{1}{\sqrt{T}})$ convergence rate result can also be achieved with a corresponding constant $C_1'=\f{\sqrt{2}d_{\phi}^2}{\eta}+[\f{8\sqrt{2}N^2(n+4)^2\Gamma\widehat{L}^2\eta}{(1-\gamma)\sigma_{\phi}}+\f{\sqrt{2}(16N+1)(n+4)^2\widehat{L}^2\eta}{\sigma_{\phi}
 }]$. The mathematical procedure is similar with the Theorem \ref{main} and Theorem \ref{crate1}. Choice on zeroth-order gradient estimators depends on the underlying science or engineering background. That is, when the information of the communication output $y_i^t$ is easy to obtain, $g_{\mu_i}(y_i^t)$ can be used; when the local information of agents $x_i^t$ is available easily, then $g_{\mu_i}(x_i^t)$ can also be used.
\end{rmk}

Till now,  the distributed randomized gradient-free mirror descent algorithm and a corresponding convergence rate are established. The algorithm has generalized some earlier work in different aspects. This section is concluded by listing the following result as a direct corollary when distance generating function is chosen by $\phi(x)=\f{1}{2}\|x\|^2$. In this case, Bregman divergence becomes $D_{\phi}(x,y)=\f{1}{2}\|x-y\|^2$, and the Bregman projection degenerates to the classical Euclidean projection.
\begin{cor} Under Assumption \ref{as2}, let $\{x_i^t\}_{t\geq 0}$ and $\{y_i^t\}_{t\geq 0}$ be the sequences generated by following projection algorithm,
\begin{eqnarray}
y_i^t&=&\sum_{j=1}^N[P^t]_{ij}x_j^t, \label{gg1} \\
 x_i^{t+1}&=&P_{X}(y_i^t-\aaa_t\f{f_i(y_i^t+\mu_i\xi_i^t)-f_i(y_i^t)}{\mu_i}\xi_i^t).\label{gg2}
\end{eqnarray}
Take stepsize by $\aaa_t=\f{\eta}{\sqrt{t+1}}$. Then for any optimal point $x^{\star}$, $\widehat{x}_l^T$ converges asymptotically to the approximate optimal solution with convergence rate $O(\f{1}{\sqrt{T}})$ with the optimal value error less than $\sqrt{n}\widehat{L}\cdot\f{1}{N}\sum_{i=1}^N\mu_i$.
\end{cor}
\begin{rmk}
Algorithm (\eqref{gg1})-(\eqref{gg2}) is exactly the Euclidean gradient-free projection algorithm in \cite{gf}. By investigating the DRGFMD algorithm,  the Euclidean gradient-free projection algorithm in \cite{gf} has already been extended to a more general Bregman non-Euclidean framework via an essentially different approach.
\end{rmk}

\subsection{DRGFMD algorithm for strongly convex optimization}
The distributed strongly-convex optimization problem is investigated by using DRGFMD in this section. Here, $X$ denotes a closed convex domain (not necessarily compact) in this section. First of all, a basic assumption on strongly convexity  is given below.
\begin{ass}\label{asst}
For each $i\in V$, $f_i: \mathbb R^n\to \mathbb R$ is assumed to be $\sigma_f$-strongly convex.
\end{ass}
The $\sigma_f$-strong convexity of $f_{\mu_i}$ is ensured in the following lemma.
\begin{lem}\label{stfmu}
Let Assumption \ref{asst} hold. Then, $f_{\mu_i}: \mathbb R^n\to \mathbb R$ is $\sigma_f$-strongly convex.
\end{lem}
\pf For any $x,y\in \mathbb R^n$ and any $\ttt\in [0,1]$, use the definition of $f_{\mu_i}$, it follows that
\begin{eqnarray}
 \nonumber &&f_{\mu_i}(\ttt x+(1-\ttt)y)   \\
 \nonumber &=&\f{1}{\kappa}\int_{\mathbb{R}^n}f_i[\ttt x+(1-\ttt)y+\mu_i\xi]e^{-\f{1}{2}\|\xi\|^2}d\xi,\\
\nonumber &=&\f{1}{\kappa}\int_{\mathbb{R}^n}f_i[\ttt (x+\mu_i\xi)+(1-\ttt)(y+\mu_i\xi)]e^{-\f{1}{2}\|\xi\|^2}d\xi, \\
\nonumber&\leq&\f{1}{\kappa}\int_{\mathbb{R}^n}[\ttt f_i(x+\mu_i\xi)+(1-\ttt)f_i(y+\mu_i\xi)\\
\nono&&-\f{\sigma_f\ttt(1-\ttt)}{2}\|(x+\mu_i\xi)-(y+\mu_i\xi)\|^2]e^{-\f{1}{2}\|\xi\|^2}d\xi\\
&=&\ttt f_{\mu_i}(x)+(1-\ttt)f_{\mu_i}(y)-\f{\sigma_f\ttt(1-\ttt)}{2}\|x-y\|^2,
\end{eqnarray}
in which the inequality follows from the convexity of $f_i$ and the third equality follows from $\f{1}{\kappa}\int_{\mathbb{R}^n}e^{-\f{1}{2}\|\xi\|^2}d\xi=1$. The proof is completed.

On the other hand, in this section, an assumption on distance generating function is given below to handle strongly convex problem.
\begin{ass}\label{distance}
Let Assumption \ref{as1} hold, the Bregman distance generating function $\phi$ is assumed to have Lipschitz gradient on $X$ with constant $\ww{L}_{\phi}$, $i. e. $
\begin{eqnarray}
\|\nabla \phi(x)-\nabla \phi(y)\|\leq \ww{L}_{\phi}\|x-y\|,\ for \ all \ x, y\in X.
\end{eqnarray}
\end{ass}
The assumption is out of consideration for practical application and following theoretical calculation. In several usual applications like machine learning, a distance generating function $\phi$ (such as $\phi(x)=\f{1}{2}\|x\|^2$ on $\mathbb R^n$ and $\phi(x)=\sum_{d=1}^n[x]_d\ln[x]_d$ on given bounded domain) can always be chosen such that $\nabla{\phi}$ is Lipschitz. The proposed strongly convex results are suitable for these cases.  Under Assumption \ref{distance}, a basic lemma for this section holds as follow,
\begin{lem} \label{DDD}
Let Assumption \ref{distance} hold, then the Bregman divergence satisfies the following relation,
\begin{eqnarray}
\nonumber D_{\phi}(x,y)\leq \f{\ww{L}_{\phi}}{2}\|x-y\|^2 \ for \ all \ x, y\in X.
\end{eqnarray}
\end{lem}
\pf Start from the definition of $D_{\phi}(x, y)$,
\begin{eqnarray}
\nonumber &&D_{\phi}(x,y)\\
\nono&=& \int_0^1\nn\nabla\phi(t(x-y)+y), x-y\mm dt-\nn\nabla\phi(y), x-y\mm   \\
\nonumber&\leq&\int_0^1\|\nabla\phi(t(x-y)+y)-\nabla\phi(y)\|\cdot\|x-y\|dt  \\
\nonumber&\leq&\int_0^1\ww{L}_{\phi}t\|x-y\|^2dt=\f{\ww{L}_{\phi}}{2}\|x-y\|^2,
\end{eqnarray}
in which the first inequality follows from Cauchy inequality, the second inequality follows from Assumption \ref{distance}.

Now it's ready to give the strongly convex convergence result for this section. For convenience of several calculations, this part uses the DRGFMD algorithm with gradient estimator $g_{\mu}(x_i^t)$.
\begin{thm}\label{strongconv}
Let Assumptions \ref{as1}, \ref{as2}, \ref{asst}, \ref{distance} hold, let $\{x_i^t\}_{t\geq 0}$ and $\{y_i^t\}_{t\geq 0}$ be the sequences generated by the DRGFMD algorithm. Let  stepsize $\aaa_0=\f{\ww{L}_{\phi}}{\sigma_f}$ and $\aaa_t=\f{\ww{L}_{\phi}}{\sigma_ft}$ for $t\geq1$, $\ww{L}_{\phi}$ is the Lipschitz constant in Assumption \ref{distance}. Then for any $l\in V$ and any optimal point $x^{\star}$,  the algorithm achieves an $O(\f{\ln T}{T})$ approximate convergence rate in following sense:
\begin{eqnarray}
\nonumber  \mathbb E[f(\widehat{x}_l^T)]-f(x^{\star})\leq \sqrt{n}\widehat{L}\cdot\f{1}{N}\sum_{i=1}^N\mu_i+C_2\cdot \f{\ln T}{T}, \ T\geq 8,
\end{eqnarray}
in which
\begin{eqnarray}
\nonumber C_2=[\f{4N^2(n+4)^2\Gamma\widehat{L}^2}{(1-\gamma)\sigma_{\phi} }+\f{(16N+1)(n+4)^2\widehat{L}^2}{2\sigma_{\phi} }]\cdot\f{2\ww{L}_{\phi}}{\sigma_f}.
\end{eqnarray}
\end{thm}

\pf Start with the inner product estimate,
\begin{small}
\begin{eqnarray}
\nonumber &&\f{1}{NT}\sum_{t=1}^T\sum_{i=1}^N\mathbb E[\nn\nabla f_{\mu_i}(x_i^t),y_i^t-x^{\star}\mm]  \\
\nonumber &=& \f{1}{NT}\sum_{t=1}^T\sum_{i=1}^N\big(\mathbb E[\nn\nabla f_{\mu_i}(x_i^t),y_i^t-x_i^t\mm]\\
\nono&&+\mathbb E[\nn\nabla f_{\mu_i}(x_i^t),x_i^t-x^{\star}\mm]\big) \\
\nonumber &\geq& \f{1}{NT}\sum_{t=1}^T\sum_{i=1}^N\mathbb E[\nn\nabla f_{\mu_i}(x_i^t),y_i^t-x_i^t\mm]\\
\nono&&+\f{1}{NT}\sum_{t=1}^T\sum_{i=1}^N\mathbb E[f_{\mu_i}(x_i^t)-f_{\mu_i}(x^{\star})] \\
\nono&&+\f{\sigma_{f}}{\ww{L}_{\phi}NT}\sum_{t=1}^T\sum_{i=1}^N\mathbb E[D_{\phi}(x^{\star},x_i^t)]\\
&&=h_1+h_2+h_3, \label{hhh}
\end{eqnarray}
\end{small}
in which the inequality follows from Lemma \ref{stfmu} and Lemma \ref{DDD}. Now the estimate for $h_1$ is given as follow,
\begin{small}
\begin{eqnarray}
\nonumber &&h_1\geq-\f{(n+4)\widehat{L}}{NT}\sum_{t=1}^T\sum_{i=1}^N\mathbb E[\|y_i^t-x_i^t\|] \\
\nonumber&&\geq -\f{(n+4)\widehat{L}}{NT}\sum_{i=1}^N\sum_{t=1}^T\sum_{j=1}^N\mathbb E[\|x_i^t-x_j^t\|]  \\
&&\geq -[\f{2N^2\Gamma(n+4)^2\widehat{L}^2 }{\sigma_{\phi}(1-\gamma)}+\f{4N(n+4)^2\widehat{L}^2}{\sigma_{\phi}}]\f{1}{T}\sum_{t=0}^T\aaa_t,  \label{h1}
\end{eqnarray}
\end{small}

\noindent in which the first inequality follows from (\eqref{kexi}), the second inequality follows from the fact that $\|\sum_{j=1}^N[P^t]_{ij}x_j^t-x_i^t\|\leq\sum_{j=1}^N[P^t]_{ij}\|x_i^t-x_j^t\|\leq\sum_{j=1}^N\|x_i^t-x_j^t\|$, and the third inequality follows from Lemma \ref{consensus}. On the other hand, for any $l\in V$,
\begin{small}
\begin{eqnarray}
\nonumber h_2&=&\f{1}{NT}\sum_{t=1}^T\sum_{i=1}^N\mathbb E[f_{\mu_i}(x_l^t)+f_{\mu_i}(x_i^t)-f_{\mu_i}(x_l^t)]-f_{\mu}(x^{\star})  \\
\nonumber&\geq&\f{1}{T}\sum_{t=1}^T\mathbb E[f_{\mu}(x_l^t)]-\f{(n+4)\widehat{L}}{NT}\sum_{t=1}^T\sum_{i=1}^N\mathbb E[\|x_i^t-x_l^t\|]-f_{\mu}(x^{\star}) \\
\nonumber&\geq&\f{1}{T}\sum_{t=1}^T \mathbb E[f(x_l^t)]-\f{(n+4)\widehat{L}}{T}\sum_{t=1}^T\sum_{i=1}^N\mathbb E[\|x_i^t-x_l^t\|]\\
\nono&&-f(x^{\star})-\f{1}{N}\sum_{i=1}^N\mu_i   \\
\nono&\geq& \mathbb E[f(\widehat{x}_l^T)]-f(x^{\star})-\f{1}{N}\sum_{i=1}^N\mu_i\\
&&-[\f{2N^2\Gamma(n+4)^2\widehat{L}^2 }{\sigma_{\phi}(1-\gamma)}+\f{4N(n+4)^2\widehat{L}^2}{\sigma_{\phi}}]\cdot\f{1}{T}\sum_{t=0}^T\aaa_t,   \label{h2}
\end{eqnarray}
\end{small}
\noindent in which the first inequality follows from the definition of $f_{\mu}$ and (\eqref{kexi}), the second inequality follows from (a) in Lemma \ref{lip} and $\f{1}{NT}\leq \f{1}{T}$, the third inequality is as a result of the convexity of $f$ and Lemma \ref{consensus}. Now combine (\eqref{hhh}), (\eqref{h1}), (\eqref{h2}) and Theorem \ref{neiji}, it can be obtained that
\begin{small}
\begin{eqnarray}
\nono&&\mathbb E[f(\widehat{x}_l^T)]-f(x^{\star})\leq \sqrt{n}\widehat{L}\cdot\f{1}{N}\sum_{i=1}^N\mu_i \\
\nono&&+\f{1}{NT}\sum_{i=1}^N\big[\f{1}{\aaa_1}\mathbb E [D_{\phi}(x^{\star},x_i^1)]   \\
\nonumber &&+\sum_{t=2}^T\mathbb E[D_{\phi}(x^{\star}, x_i^t)](\f{1}{\aaa_t}-\f{1}{\aaa_{t-1}})-\f{1}{\aaa_T}\mathbb E[D_{\phi}(x^{\star}, x_i^{T+1})]\big]\\
\nono&&-\f{\sigma_{f}}{\ww{L}_{\phi}NT}\sum_{t=1}^T\sum_{i=1}^N\mathbb E[D_{\phi}(x^{\star},x_i^t)]+\\
\nono&&[\f{4N^2\Gamma(n+4)^2\widehat{L}^2 }{\sigma_{\phi}(1-\gamma)}+\f{(16N+1)(n+4)^2\widehat{L}^2}{2\sigma_{\phi}}]\cdot\f{1}{T}\sum_{t=0}^T\aaa_t\\
&=&l_1+l_2-l_3+l_4.
\end{eqnarray}
\end{small}
Since $\aaa_t=\f{\ww{L}_{\phi}}{\sigma_{f}t}$ when $t\geq1$ and the Bregman divergence is non-negative, it follows that
\begin{eqnarray}
\nonumber l_2-l_3&=& \f{1}{NT}(\f{1}{\aaa_1}-\f{\sigma_f}{\ww{L}_{\phi}})\sum_{i=1}^N\mathbb E[D_{\phi}(x^{\star},x_i^1)]   \\
\nonumber&&+\f{1}{NT}\sum_{t=2}^T(\f{1}{\aaa_t}-\f{1}{\aaa_{t-1}}-\f{\sigma_{f}}{\ww{L}_{\phi}})\sum_{i=1}^N\mathbb E[D_{\phi}(x^{\star},x_i^t)] \\
\nono&&-\f{1}{NT\aaa_T}\sum_{i=1}^N\mathbb E[D_{\phi}(x^{\star},x_i^{T+1})]\\
&=&-\f{1}{NT\aaa_T}\sum_{i=1}^N\mathbb E[D_{\phi}(x^{\star},x_i^{T+1})]\leq 0.
\end{eqnarray}
Therefore,
\begin{small}
\begin{eqnarray}
 \nono&&\mathbb E[f(\widehat{x_i}^T)]-f(x^{\star})\leq \sqrt{n}\widehat{L}\cdot\f{1}{N}\sum_{i=1}^N\mu_i\\
\nono &&+[\f{4N^2(n+4)^2\Gamma\widehat{L}^2}{(1-\gamma)\sigma_{\phi} }+\f{(16N+1)(n+4)^2\widehat{L}^2}{2\sigma_{\phi} }]\cdot\f{1}{T}\sum_{t=0}^T\aaa_t.
\end{eqnarray}
\end{small}

\noindent In addition, use the fact that $\aaa_0=\f{\ww{L}_{\phi}}{\sigma_f}$ and $\aaa_t=\f{\ww{L}_{\phi}}{\sigma_f t}$ for $t\geq1$, it follows that
\begin{eqnarray}
\nono&&\f{1}{T}\sum_{t=0}^T\aaa_t=\f{1}{T}(\f{\ww{L}_{\phi}}{\sigma_f}+\sum_{t=1}^T\f{\ww{L}_{\phi}}{\sigma_ft})\\
&&\leq \f{\ww{L}_{\phi}}{\sigma_f T}(2+\int_{1}^T\f{1}{s}ds)\leq \f{2\ww{L}_{\phi}}{\sigma_f}\cdot \f{\ln T}{T}, T\geq 8,     \label{stest}
\end{eqnarray}
then the desired result holds.

\begin{rmk}
Theorem \ref{strongconv} shows that the DRGFMD algorithm achieves an $O(\f{\ln T}{T})$ approximate convergence rate for strongly convex constrained optimization over time-varying network, generalizing the one in \cite{strongc} to a non-Euclidean situation. In addition, any smoothness assumptions on objective functions are not needed. Moreover, the proposed algorithm is the first distributed non-Euclidean zeroth-order method applied to strongly convex optimization problem and the smoothing function $f_{\mu_i}$ acts as an important bridge to achieve the final convergence rate.
\end{rmk}

\section{The DRGFMD algorithm with weighted averaging}
The former sections of the paper have discussed the DRGFMD algorithm in distributed convex and strongly optimization problem. However, the approximating sequence to the convergence of algorithm in former sections are all in classical form $\widehat{x}_l^T=\f{1}{T}\sum_{t=1}^T x_l^{t}$. The paper in this section provides a weighted average approximating sequence which is different from the approximating sequence the existing distributed algorithms have used. In this section, the DRGFMD algorithm with weighted averaging (\textbf{DRGFMD-WA}) is investigated and applied to the convex and strongly convex optimization problem. Several estimates obtained in the former section will be used directly in this section. For convenience of saving space and without loss of generality, $x_i^0=0$ for all $i\in V$ is still assumed in this section. By setting $\Phi_v^1=\max_{i,j\in V}\|x_i^0-x_j^0\|$ and $\Phi_v^2=\max_{i\in V}\|x_i^0\|$ in corresponding constant place, convergence results for non-zero initial data case can be gotten.
\subsection{DRGFMD compact constrained convex optimization with weighted averaging}
For any $l\in V$, denote the weighted average approximating sequence by
\begin{eqnarray}
 \widetilde{x}_l^T=\f{\sum_{t=1}^T\f{x_l^{t}}{\aaa_t} }{\sum_{t=1}^T\f{1}{\aaa_t}}, \label{weightseq}
\end{eqnarray}
then the first distributed convergence result with weighted average approximating sequence is given as follow.
\begin{thm}\label{weight1}
Let Assumptions \ref{as1} and \ref{as2} hold. Let the stepsize $\aaa_t$ is a non-increasing sequence. Then, for the weighted average sequence $ \widetilde{x}_l^T$  generated by the DRGFMD algorithm, the following convergence result holds for any $l\in V$ and $T\geq 1$,
\begin{eqnarray}
\nonumber &&\mathbb E[f(\widetilde{x}_l^T)]-f(x^{\star})\leq \sqrt{n}\widehat{L}\cdot\f{1}{N}\sum_{i=1}^N\mu_i\\
\nono&&+\f{1}{\sum_{t=1}^T\f{1}{\aaa_t}}[\f{d_{\phi}^2}{\aaa_T^2}+\f{(n+4)^2\widehat{L}^2T}{2\sigma_{\phi}}
+\f{\widetilde{C}}{\aaa_T}\sum_{t=0}^T\aaa_t],
\end{eqnarray}
in which
\begin{eqnarray}
\nonumber \widetilde{C}=\f{2N(n+4)^2\Gamma\widehat{L}^2}{\sigma_{\phi} (1-\gamma)}+\f{4(n+4)^2\widehat{L}^2}{\sigma_{\phi} }.
\end{eqnarray}

\end{thm}

\pf Start from (\eqref{sides}) in Theorem \ref{neiji} with $z_i^t=y_i^t$ as follow,
\begin{eqnarray}
 \nono&&\aaa_t\nn g_{\mu_i}(y_i^t), y_i^t-x^{\star}\mm\\
 &&\leq D_{\phi}(x^{\star},y_i^t)-D_{\phi}(x^{\star},x_i^{t+1})+\f{\aaa_t^2}{2\sigma_{\phi}}\|g_{\mu_i}(y_i^t)\|^2.   \label{w1}
\end{eqnarray}
Take the conditional expectation on $F_t$ on both sides of (\eqref{w1}), note the fact that $D(x^{\star},y_i^t)$ is measurable with respect to $F_t$, rearrange terms and use Lemma \ref{lip}, it follows that
\begin{eqnarray}
 \nono&& \aaa_t\nn \nabla f_{\mu_i}(y_i^t), y_i^t-x^{\star}\mm+\mathbb E[D_{\phi}(x^{\star},x_i^{t+1})| F_t]\\
&& \leq D_{\phi}(x^{\star},y_i^t)+\f{\aaa_t^2(n+4)^2\widehat{L}^2}{2\sigma_{\phi}}.    \label{w2}
\end{eqnarray}
Dividing both sides of the above inequality by $\aaa_t^2$ and rearranging terms yields
\begin{eqnarray}
 \nono &&\f{1}{\aaa_t}\nn  \nabla f_{\mu_i}(y_i^t), y_i^t-x^{\star}\mm+\f{1}{\aaa_t^2}\mathbb E[D_{\phi}(x^{\star},x_i^{t+1})| F_t]\\
 &&\leq \f{1}{\aaa_t^2}D_{\phi}(x^{\star},y_i^t)+\f{(n+4)^2\widehat{L}^2}{2\sigma_{\phi}}. \label{div}
\end{eqnarray}
Note that
\begin{eqnarray}
\nonumber\f{1}{\aaa_t^2}D_{\phi}(x^{\star},y_i^t)&=& \f{1}{\aaa_{t-1}^2}D_{\phi}(x^{\star},y_i^t)+(\f{1}{\aaa_t^2}-\f{1}{\aaa_{t-1}^2})D_{\phi}(x^{\star},y_i^t) \\
&\leq& \f{1}{\aaa_{t-1}^2}D_{\phi}(x^{\star},y_i^t)+(\f{1}{\aaa_t^2}-\f{1}{\aaa_{t-1}^2})d_{\phi}^2, \label{div1}
\end{eqnarray}
substitute (\eqref{div1}) into (\eqref{div}) and take total expectation on both sides, it follows that
\begin{eqnarray}
\nono&& \f{1}{\aaa_t} \mathbb E[\nn \nabla f_{\mu_i}(y_i^t), y_i^t-x^{\star}\mm]+\f{1}{\aaa_t^2}\mathbb E[D_{\phi}(x^{\star},x_i^{t+1})]\\
 \nono&&\leq \f{1}{\aaa_{t-1}^2}\mathbb E[D_{\phi}(x^{\star},y_i^t)]+(\f{1}{\aaa_t^2}-\f{1}{\aaa_{t-1}^2})d_{\phi}^2+\f{(n+4)^2\widehat{L}^2}{2\sigma_{\phi}}.
\end{eqnarray}
Set $\aaa_0=1$ and sum up both sides from $t=1$ to $T$, it follows that
\begin{eqnarray}
 \nonumber &&\sum_{t=1}^T\f{1}{\aaa_t} \mathbb E[\nn \nabla f_{\mu_i}(y_i^t), y_i^t-x^{\star}\mm]+\sum_{t=1}^T\f{1}{\aaa_t^2}\mathbb E[D_{\phi}(x^{\star},x_i^{t+1})]\\
\nono &&\leq \sum_{t=1}^T\f{1}{\aaa_{t-1}^2}\mathbb E[D_{\phi}(x^{\star},y_i^t)]+(\f{1}{\aaa_T^2}-1)d_{\phi}^2+\f{(n+4)^2\widehat{L}^2T}{2\sigma_{\phi}}.
\end{eqnarray}
Sum up both sides from $i=1$ to $N$ and substitute $y_i^t=\sum_{j=1}^N[P^t]_{ij}x_j^t$ into the right hand side, it can be obtained that
\begin{small}
\begin{eqnarray}
\nono&& \sum_{t=1}^T\f{1}{\aaa_t} \sum_{i=1}^N\mathbb E[\nn \nabla f_{\mu_i}(y_i^t), y_i^t-x^{\star}\mm]  \\
\nono&&+\sum_{i=1}^N\sum_{t=1}^T\f{1}{\aaa_t^2}\mathbb E[D_{\phi}(x^{\star},x_i^{t+1})]\leq\sum_{i=1}^N\sum_{t=1}^T\f{1}{\aaa_{t-1}^2}\mathbb E[D_{\phi}(x^{\star},x_i^t)]\\
\nono&&+N(\f{1}{\aaa_T^2}-1)d_{\phi}^2+\f{N(n+4)^2\widehat{L}^2T}{2\sigma_{\phi}},
\end{eqnarray}
\end{small}
in which the inequality follows from the separate convexity of Bregman divergence $D_{\phi}(x,y)$ and the doubly stochastic property of the communication matrix $P^t$. Now delete the same terms of both sides, it follows that
\begin{small}
\begin{eqnarray}
&& \nonumber \sum_{t=1}^T\f{1}{\aaa_t} \sum_{i=1}^N\mathbb E[\nn \nabla f_{\mu_i}(y_i^t), y_i^t-x^{\star}\mm]+\sum_{i=1}^N\f{\mathbb E[D_{\phi}(x^{\star},x_i^{T+1})]}{\aaa_T^2}  \\
\nonumber &\leq&\sum_{i=1}^N\mathbb E[D_{\phi}(x^{\star},x_i^1)]+N(\f{1}{\aaa_T^2}-1)d_{\phi}^2+\f{N(n+4)^2\widehat{L}^2T}{2\sigma_{\phi}}  \\
\nonumber&\leq& Nd_{\phi}^2+N(\f{1}{\aaa_T^2}-1)d_{\phi}^2+\f{N(n+4)^2\widehat{L}^2T}{2\sigma_{\phi}}\\
&=&\f{Nd_{\phi}^2}{\aaa_T^2}+\f{N(n+4)^2\widehat{L}^2T}{2\sigma_{\phi}},
\end{eqnarray}
\end{small}
in which the second inequality follows from the compactness of $X$. Since $\sum_{i=1}^N\f{1}{\aaa_T^2}\mathbb E[D_{\phi}(x^{\star},x_i^{T+1})]$ is nonnegative, this fact leads to
\begin{small}
\begin{eqnarray}
  \sum_{t=1}^T\f{1}{\aaa_t} (\f{1}{N}\sum_{i=1}^N\mathbb E[\nn \nabla f_{\mu_i}(y_i^t), y_i^t-x^{\star}\mm])\leq \f{d_{\phi}^2}{\aaa_T^2}+\f{(n+4)^2\widehat{L}^2T}{2\sigma_{\phi}}.   \label{neicon}
\end{eqnarray}
\end{small}
Combine (\eqref{neicon}) and (\eqref{fmu2}), it can be obtained that, for any $l\in V$,
\begin{eqnarray}
\nonumber &&\sum_{t=1}^T\f{1}{\aaa_t}(\mathbb E[f(x_l^t)]-f(x^{\star}))\\
\nono&&\leq \sum_{t=1}^T\f{1}{\aaa_t}\sqrt{n}\widehat{L}\cdot\f{1}{N}\sum_{i=1}^N\mu_i+\f{d_{\phi}^2}{\aaa_T^2}+\f{(n+4)^2\widehat{L}^2T}{2\sigma_{\phi}}\\
&&+\f{(n+4)\widehat{L}}{N}\sum_{t=1}^T\f{1}{\aaa_t}\cdot\sum_{j=1}^N\mathbb E[\|x_j^t-x_l^t\|]. \label{conccc}
\end{eqnarray}
use the non-increasing assumption of $\aaa_t$ and the consensus result Lemma \ref{consensus}, the following holds,
\begin{eqnarray}
\nonumber&&\f{(n+4)\widehat{L}}{N}\sum_{t=1}^T\f{1}{\aaa_t}\cdot\sum_{j=1}^N\mathbb E[\|x_j^{t}-x_l^{t}\|]\\
\nonumber&\leq& \f{(n+4)\widehat{L}}{N}\f{1}{\aaa_T}\cdot\sum_{t=1}^T\sum_{j=1}^N\mathbb E[\|x_j^{t}-x_l^{t}\|] \\
&\leq&(\f{2N(n+4)^2\Gamma\widehat{L}^2}{\sigma_{\phi} (1-\gamma)}+\f{4(n+4)^2\widehat{L}^2}{\sigma_{\phi} })\f{1}{\aaa_T}\sum_{t=0}^T\aaa_t.\label{conccc1}
\end{eqnarray}
The desired result follows by combining (\eqref{conccc}) with (\eqref{conccc1}), dividing by $\sum_{t=1}^T\f{1}{\aaa_t}$ on both sides and noticing the convexity of $f$.
\begin{rmk}
In addition to the smoothing parameter term, the convergence result consists of three terms, they are all under the influence of the weighted averaging. The first term $\f{d_{\phi}^2}{\aaa_T^2}$ represents a topology effect from the underlying space $X$, the second term $\f{(n+4)^2\widehat{L}^2T}{2\sigma_{\phi}}$ represents an intrinsic centralized effect, the third term $\f{\widetilde{C}}{\aaa_T}\sum_{t=0}^T\aaa_t$ is the decentralized term which is as a result of the network topology and the distributed communication of information in the network.
\end{rmk}

In what follows, the convergence rate is considered. Let $0<\delta<1$, the stepsize of the following form is used:
\begin{eqnarray}
\aaa_t=\f{\rho}{(t+1)^{\delta}},t\geq1 \ and \ \aaa_0=1. \label{step}
\end{eqnarray}
Before obtaining the convergence rate, the following inequality of $\sum_{t=1}^T\f{1}{\aaa_t}$ is needed for providing a lower bound estimate.
\begin{lem}\label{xulie}
Let the stepsize $\aaa_t$ be defined as (\eqref{step}), for any $T\geq1$ and $p\leq 1-\f{1}{2^{1+\delta}}$, the following estimate holds,
\begin{eqnarray}
\sum_{t=1}^T\f{1}{\aaa_t}\geq \f{p(T+1)^{\delta+1}}{\rho(\delta+1)}.
\end{eqnarray}
\end{lem}
\pf
According to the concavity of function $s(t)=(t+1)^{\delta}$ for $t\geq1$, the following holds,
\begin{eqnarray}
\nono&&\sum_{t=1}^T\f{1}{\aaa_t}=\sum_{t=1}^T\f{1}{\rho}(t+1)^{\delta}\geq\f{1}{\rho}\int_{0}^{T}(t+1)^{\delta}dt\\
&&=\f{1}{\rho(\delta+1)}((T+1)^{\delta+1}-1).
\end{eqnarray}
Select a $p$ such that $p\leq1-\f{1}{2^{\delta+1}}$, then $(T+1)^{\delta+1}-1\geq p(T+1)^{\delta+1}$ for any $T\geq1$ and the desired result holds.

With convergence result Theorem \ref{weight1} and lower estimate Lemma \ref{xulie} in hand, it's ready to present the convergence rate result.

\begin{thm}\label{weightrateder}
Let Assumptions \ref{as1} and \ref{as2} hold. Let the stepsize $\aaa_t$ be the sequence given by (\eqref{step}). Then for the weighted average sequence $ \widetilde{x}_l^T$ generated by the DRGFMD algorithm and all $T\geq 1$,
\begin{eqnarray}
\nonumber \mathbb E[f(\widetilde{x}_l^T)]-f(x^{\star})&\leq& \sqrt{n}\widehat{L}\cdot\f{1}{N}\sum_{i=1}^N\mu_i+C_{\delta, 1}\cdot\f{1}{(T+1)^{1-\delta}}\\
&&+C_{\delta, 2}\cdot\f{1}{(T+1)^{\delta}}.  \label{weightcon}
\end{eqnarray}
in which
\begin{eqnarray}
\nonumber C_{\delta, 1}&=&\f{(\delta+1)d_{\phi}^2}{p\rho},\\
\nonumber C_{\delta, 2}&=&[\f{(n+4)^2\widehat{L}^2}{2\sigma_{\phi}}+\f{\widetilde{C}(\rho+1-\delta)}{\rho(1-\delta)}]\cdot\f{\rho(\delta+1)}{p},\\
\nonumber \widetilde{C}&=&\f{2N(n+4)^2\Gamma\widehat{L}^2}{\sigma_{\phi} (1-\gamma)}+\f{4(n+4)^2\widehat{L}^2}{\sigma_{\phi} }.
\end{eqnarray}
\end{thm}
\pf Combine Theorem \ref{weight1} and Lemma \ref{xulie}, and notice that the following fact holds,
\begin{eqnarray}
\nono&&\sum_{t=0}^T\aaa_t\leq1+\int_{0}^T\f{\rho}{(1+t)^{\delta}}\leq \f{\rho(1+T)^{1-\delta}}{1-\delta}-\f{\rho}{1-\delta}+1\\
\nono&&\leq\f{(\rho+1-\delta)(1+T)^{1-\delta}}{1-\delta},
\end{eqnarray}
then it follows that
\begin{eqnarray}
\nonumber &&\f{1}{\sum_{t=1}^T\f{1}{\aaa_t}}[\f{d_{\phi}^2}{\aaa_T^2}+\f{(n+4)^2\widehat{L}^2T}{2\sigma_{\phi}}
+\f{\widetilde{C}}{\aaa_T}\sum_{t=0}^T\aaa_t]\\
\nonumber&\leq& \f{\rho(\delta+1)}{p(T+1)^{\delta+1}}[\f{d_{\phi}^2(T+1)^{2\delta}}{\rho^2}+\f{(n+4)^2\widehat{L}^2(T+1)}{2\sigma_{\phi}}\\
\nono&&+\f{\widetilde{C}(\rho+1-\delta)(T+1)}{\rho(1-\delta)}]\\
\nonumber&=& \f{(\delta+1)d_{\phi}^2}{p\rho}\cdot\f{1}{(T+1)^{1-\delta}}\\
\nono&&+ [\f{(n+4)^2\widehat{L}^2}{2\sigma_{\phi}}+\f{\widetilde{C}(\rho+1-\delta)}{\rho(1-\delta)}]\cdot\f{\rho(\delta+1)}{p}\cdot\f{1}{(T+1)^{\delta}},
\end{eqnarray}
which implies (\eqref{weightcon}) and the proof is concluded.
\begin{cor}\label{DDDD}
Under assumptions of Theorem \ref{weightrateder}, let $C_{\delta}=\max\{C_{\delta, 1}, C_{\delta, 2}\}$, in which $C_{\delta, 1}$ and $C_{\delta, 2}$ are the constants in Theorem \ref{weightrateder}. For the weighted average sequence $ \widetilde{x}_l^T$ generated by the DRGFMD algorithm and  all $T\geq 1$, the following approximate convergence rate for the DRGFMD algorithm holds:
\begin{small}
\begin{eqnarray}
\nonumber  \mathbb E[f(\widetilde{x}_l^T)]-f(x^{\star})&\leq& \sqrt{n}\widehat{L}\cdot\f{1}{N}\sum_{i=1}^N\mu_i+\f{2C_{\delta}}{T^{\delta}} \ if \ \delta\in (0,\f{1}{2});\\
\nonumber  \mathbb E[f(\widetilde{x}_l^T)]-f(x^{\star})&\leq& \sqrt{n}\widehat{L}\cdot\f{1}{N}\sum_{i=1}^N\mu_i+\f{2C_{\delta}}{\sqrt{T}} \ if \ \delta=\f{1}{2};\\
\nonumber  \mathbb E[f(\widetilde{x}_l^T)]-f(x^{\star})&\leq& \sqrt{n}\widehat{L}\cdot\f{1}{N}\sum_{i=1}^N\mu_i+\f{2C_{\delta}}{T^{1-\delta}} if \ \delta\in (\f{1}{2},1).
\end{eqnarray}
\end{small}
\end{cor}
\pf Note that
\begin{small}
\begin{eqnarray}
\nonumber \f{1}{(T+1)^{1-\delta}}\leq\f{1}{(T+1)^{\delta}}\leq \f{1}{T^{\delta}} \ if  \ \delta \in (0,\f{1}{2}); \\
\nonumber \f{1}{(T+1)^{\delta}}=\f{1}{(T+1)^{1-\delta}}=\f{1}{(T+1)^{\f{1}{2}}}\leq \f{1}{\sqrt{T}} \ \ if \ \delta=\f{1}{2}; \\
\nono \f{1}{(T+1)^{\delta}}\leq\f{1}{(T+1)^{1-\delta}}\leq \f{1}{T^{1-\delta}}  \ if \ \delta\in (\f{1}{2},1),
\end{eqnarray}
\end{small}
then by taking $C_{\delta}=\max\{C_{\delta, 1}, C_{\delta, 2}\}$, the desired results follow directly.
\begin{rmk}
Theorem \ref{weight1}, Theorem \ref{weightrateder} and Corollary \ref{DDDD} provide a general analyzing framework of distributed mirror descent with weighted average approximating sequence. To the best of our knowledge, this paper is the first to utilize a decentralized reciprocal weighted average approximating sequence  $\widetilde{x}_l^T=\sum_{t=1}^T\f{x_l^t}{\aaa_t}/\sum_{t=1}^T\f{1}{\aaa_t}$ to achieve a class of distributed convergence rates when stepsize is taken in $\f{\rho}{(t+1)^{\delta}}$ form, in contrast to the classical class of approximating sequence form $\widehat{x_i}^T=\f{1}{T}\sum_{t=1}^T x_i^{t}$ and $\overline{x}^t=\f{1}{N}\sum_{i=1}^Nx_i^t$ that the existing distributed mirror descent methods used. Thus the approach in this section has shed light on investigation of different types of decentralized weighted average sequences and corresponding convergence rates.
\end{rmk}

\subsection{DRGFMD strongly convex optimization with weighted averaging}
In this section, the DRGFMD algorithm with weighted averaging is used to solve the strongly convex optimization problem on the convex (not necessarily compact) constraint set. In order to construct a decentralized weighted average approximating sequence for strongly convex case, the following scaling variant of the DRGFMD algorithm (\textbf{DRGFMD$'$}) is used in this section:
\begin{small}
\begin{eqnarray}
 \nonumber y_i^t&=&\sum_{j=1}^N[P^t]_{ij}x_j^t, \\
\nonumber  x_i^{t+1}&=&\arg\min_{x\in X}\{\f{\ww{L}_{\phi}\aaa_t}{\sigma_{f}}\nn \f{f_i(y_i^t+\mu_i\xi_i^t)-f_i(y_i^t)}{\mu_i}\xi_i^t,x\mm+D_{\phi}(x,y_i^t)\},
\end{eqnarray}
\end{small}
in which $\sigma_{f}$ is the common strongly convex constant of $f_i$, $\ww{L}_{\phi}$ is the constant in Assumption \ref{distance}.

\begin{thm}\label{strw}
Let Assumptions \ref{as1}, \ref{as2}, \ref{asst}, \ref{distance} hold, let $\aaa_t$ be a non-increasing positive sequence satisfying $\aaa_0=1$ and $\f{1-\aaa_t}{\aaa_t^2}\leq\f{1}{\aaa_{t-1}^2}$ for $t\geq 1$. Let the weighted average sequence $ \widetilde{x}_l^T$  be generated by the DRGFMD$'$ algorithm. Denote $\Delta_{\phi}^0=D_{\phi}(x^{\star}, 0)$,  then for any $l\in V$ and $T\geq1$, the following strongly convex convergence result holds,
\begin{eqnarray}
 \nono&&\mathbb E[f(\widetilde{x}_l^T)]-f(x^{\star})\leq \sqrt{n}\widehat{L}\cdot\f{1}{N}\sum_{i=1}^N\mu_i\\
 \nono&&+\f{1}{\sum_{t=1}^T\f{1}{\aaa_t}}[\f{\sigma_f\Delta_{\phi}^0}{\ww{L}_{\phi}}+\f{\ww{L}_{\phi}(n+4)^2\widehat{L}^2T}{2\sigma_{f}\sigma_{\phi}}
+\f{\widetilde{C}}{\aaa_T}\sum_{t=0}^T\aaa_t],
\end{eqnarray}
where $\widetilde{C}$ is the constant in the last section.
\end{thm}
\pf Since the DRGFMD$'$ algorithm (scaling version of DRGFMD) is considered, now start from (\eqref{w2}) with $\aaa_t$ replaced by $\f{\ww{L}_{\phi}\aaa_t}{\sigma_{f}}$ and take total expectation on both sides, it follows that
\begin{eqnarray}
\nono&& \f{\ww{L}_{\phi}\aaa_t}{\sigma_{f}}\mathbb E[\nn \nabla f_{\mu_i}(y_i^t), y_i^t-x^{\star}\mm]+\mathbb E[D_{\phi}(x^{\star},x_i^{t+1})]\\
\nono &&\leq \mathbb E[D_{\phi}(x^{\star},y_i^t)]+\f{\ww{L}_{\phi}^2\aaa_t^2(n+4)^2\widehat{L}^2}{2\sigma_{f}^2\sigma_{\phi}}.
\end{eqnarray}
Sum up both sides from $i=1$ to $N$ and divide by $N$, it follows that
\begin{small}
\begin{eqnarray}
 \nono&&\f{\ww{L}_{\phi}\aaa_t}{\sigma_{f}}\f{1}{N}\sum_{i=1}^N\mathbb E[\nn \nabla f_{\mu_i}(y_i^t), y_i^t-x^{\star}\mm]\\
 \nono&&+\f{\sum_{i=1}^N\mathbb E[D_{\phi}(x^{\star},x_i^{t+1})]}{N}\\
&&\leq \f{1}{N}\sum_{i=1}^N\mathbb E[D_{\phi}(x^{\star},y_i^t)]+\f{\ww{L}_{\phi}^2\aaa_t^2(n+4)^2\widehat{L}^2}{2\sigma_{f}^2\sigma_{\phi}}. \label{222}
\end{eqnarray}
\end{small}

\noindent According to Lemma \ref{stfmu}, the strongly convex version of (\eqref{fmu2}) holds as follow,
\begin{small}
\begin{eqnarray}
\nonumber &&\f{1}{N}\sum_{i=1}^N\nn \nabla f_{\mu_i}(y_i^t), y_i^t-x^{\star}\mm \\
\nono&\geq& \f{1}{N}\sum_{i=1}^Nf_{\mu_i}(y_i^t)-\f{1}{N}\sum_{i=1}^Nf_{\mu_i}(x^{\star})+\f{1}{N}\sum_{i=1}^N\f{\sigma_{f}}{2}\|y_i^t-x^{\star}\|^2 \\
\nonumber&\geq& f(x_l^t)-f(x^{\star})-\sqrt{n}\widehat{L}\cdot\f{1}{N}\sum_{i=1}^N\mu_i-\f{(n+4)\widehat{L}}{N}\sum_{j=1}^N\|x_j^t-x_l^t\|\\
&&+\f{\sigma_{f}}{\ww{L}_{\phi}}\f{1}{N}\sum_{i=1}^ND_{\phi}(x^{\star},y_i^t),  \label{neineinei}
\end{eqnarray}
\end{small}
in which Lemma \ref{DDD} is used in the second inequality. Take total expectation of (\eqref{neineinei}) and substitute it into (\eqref{222}), after rearranging terms and dividing both sides by $\aaa_t^2$, it can be obtained that for any $l\in V$,
\begin{small}
\begin{eqnarray}
\nonumber&&\f{\ww{L}_{\phi}}{\sigma_{f}}\cdot\f{1}{\aaa_t}\mathbb E[f(x_l^t)-f(x^{\star})]+\f{1}{\aaa_t^2}\cdot\f{1}{N}\sum_{i=1}^N\mathbb E[D_{\phi}(x^{\star},x_i^{t+1})]\\
\nono&\leq& \f{1-\aaa_t}{\aaa_t^2}\f{1}{N}\sum_{i=1}^N\mathbb E[D_{\phi}(x^{\star},y_i^t)]+\f{\ww{L}_{\phi}^2(n+4)^2\widehat{L}^2}{2\sigma_{f}^2\sigma_{\phi}}\\
\nono&& +\f{\ww{L}_{\phi}}{\sigma_{f}\aaa_t}\cdot\f{(n+4)\widehat{L}}{N}\sum_{j=1}^N\mathbb E[\|x_j^t-x_l^t\|]\\
&&+\f{\ww{L}_{\phi}}{\sigma_{f}\aaa_t}\cdot\sqrt{n}\widehat{L}\cdot\f{1}{N}\sum_{i=1}^N\mu_i.      \label{dai}
\end{eqnarray}
\end{small}

\noindent Substitute $\f{1-\aaa_t}{\aaa_t^2}\leq\f{1}{\aaa_{t-1}^2}$,
$\sum_{i=1}^N\mathbb E[D_{\phi}(x^{\star},y_i^t)]\leq  \sum_{i=1}^N\sum_{j=1}^N[P^t]_{ij}\mathbb E[D_{\phi}(x^{\star},x_j^t)]=\sum_{i=1}^N\mathbb E[D_{\phi}(x^{\star},x_i^t)]$,
and consensus estimate Lemma \ref{consensus} into (\eqref{dai}), sum up both sides from $t=1$ to $T$, after using the non-increasing property of stepsize $\aaa_t$, it can be obtained that
\begin{small}
\begin{eqnarray}
\nonumber&&\f{\ww{L}_{\phi}}{\sigma_{f}}\cdot\sum_{t=1}^T\f{1}{\aaa_t}\mathbb E[f(x_l^t)-f(x^{\star})]+\f{1}{N}\sum_{i=1}^N\f{1}{\aaa_T^2}\mathbb E[D_{\phi}(x^{\star},x_i^{T+1})]\\
\nonumber&&\leq \f{1}{N}\sum_{i=1}^N\mathbb E[D_{\phi}(x^{\star},x_i^0)]+
 \f{\ww{L}_{\phi}}{\sigma_{f}}\cdot\f{\widetilde{C}}{\aaa_T}\sum_{t=0}^T\aaa_t\\
\nono&& +\f{\ww{L}_{\phi}^2(n+4)^2\widehat{L}^2T}{2\sigma_{f}^2\sigma_{\phi}}
 +\f{\ww{L}_{\phi}}{\sigma_{f}}\cdot(\sum_{t=1}^T\f{1}{\aaa_t})\cdot\f{\sqrt{n}\widehat{L}}{N}\sum_{i=1}^N\mu_i.
\end{eqnarray}
\end{small}

\noindent Note that  $\f{1}{N}\sum_{i=1}^N\f{1}{\aaa_T^2}\mathbb E[D_{\phi}(x^{\star},x_i^{T+1})]$ is nonnegative, it follows that
\begin{eqnarray}
\nonumber&&\f{\ww{L}_{\phi}}{\sigma_{f}}\cdot\sum_{t=1}^T\f{1}{\aaa_t}\mathbb E[f(x_l^t)-f(x^{\star})]
\leq \Delta_{\phi}^0+
 \f{\ww{L}_{\phi}}{\sigma_{f}}\cdot\f{\widetilde{C}}{\aaa_T}\sum_{t=0}^T\aaa_t\\
 \nono&&+\f{\ww{L}_{\phi}^2(n+4)^2\widehat{L}^2T}{2\sigma_{f}^2\sigma_{\phi}}
 +\f{\ww{L}_{\phi}}{\sigma_{f}}(\sum_{t=1}^T\f{1}{\aaa_t})\f{\sqrt{n}\widehat{L}}{N}\sum_{i=1}^N\mu_i.
\end{eqnarray}
Dividing by $\f{\ww{L}_{\phi}}{\sigma_{f}}\sum_{t=1}^T\f{1}{\aaa_t}$ on both sides of the inequality above and using the convexity of $f$ yields the desired result.

The following corollary gives an convergence rate result for the proposed distributed algorithm with weighted average approximating sequence.
\begin{cor}
Under assumptions of Theorem \ref{strw}, choose stepsize $\aaa_0=1$ and $\aaa_t=\f{2}{t+1}$ for $t\geq1$. Then $ \widetilde{x}_l^T$  generated by the DRGFMD$'$ algorithm achieves an $O(\f{\ln T}{T})$ approximate convergence rate.
\end{cor}
\pf Substitute $\aaa_t=\f{2}{t+1}$ into the right hand side of Theorem \ref{strw}, it equals to
\begin{small}
\begin{eqnarray}
\nono\f{4\sigma_{f}\Delta_{\phi}^0}{\ww{L}_{\phi}T(T+3)}+\f{2\ww{L}_{\phi}(n+4)^2\widehat{L}^2}{\sigma_{f}\sigma_{\phi}}\cdot\f{1}{T+3}+\f{2\widetilde{C}(T+1)}{T(T+3)}\sum_{t=1}^T\f{2}{t+1}.
\end{eqnarray}
\end{small}
\noindent applying the similar estimate idea of (\eqref{stest}) to the third term, the summation above equals to
\begin{eqnarray}
O(\f{1}{T^2})+O(\f{1}{T})+O(\f{\ln T}{T})=O(\f{\ln T}{T}),
\end{eqnarray}
which completes the proof.
\begin{rmk}
Till now, a strongly convex convergence rate $O(\ln T/T)$ is established for DRGFMD with decentralized weighted average approximating sequence. The procedure of canceling terms is essentially different from the existing distributed strongly convex optimization  methods, that results in the weighted average approximating sequence. Besides, both the strongly convex optimization methods in this section and last section utilize the smoothing function $f_{\mu_i}$ to serve as a bridge for proposed algorithm to convergence, which is also different from the existing strongly convex optimization methods. On the other hand, if in addition, a compact constraint assumption is permitted, Yuan et al. in \cite{md3} show that $O(\ln T/T)$ can be improved to $O(1/T)$ via an epoch distributed stochastic mirror descent method. However, there is a problem remained here: If the compact constraint condition is not satisfied, is $O(\ln T/T)$ the optimal strongly convex convergence rate? Or, does there exist a method to improve the strongly convex optimization convergence rate $O(\ln T/T)$ without the compact constraint condition? Further investigation is required to answer these questions.
\end{rmk}

\begin{rmk}
It is noteworthy that the future investigation on decentralized reciprocal weighted average sequence is necessary. As a beginning that the decentralized reciprocal weighted average sequence is applied to distributed optimization problem, it is highly possible that the decentralized reciprocal weighted average type approximation can provide some help in improving the convergence rate for other distributed optimization methods in the future.
\end{rmk}

\begin{rmk}
The convergence results in this paper are all in approximate convergence manner (up to a controllable error bound $\sqrt{n}\widehat{L}\cdot\f{1}{N}\sum_{i=1}^N\mu_i$). In fact, if the smoothing parameter $\mu_i$ is selected in a time-varying diminishing way ($\mu_i^t$), after some technical procedures, the convergence can be made to be exact convergence. It is desirable to investigate the time-varying parameters in our future work.

\end{rmk}

\section{Illustrative simulation example}
\begin{figure}[htbp]
\centering
\includegraphics[width=2in]{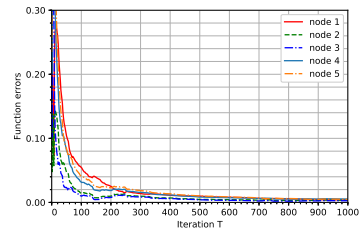}
\caption{Convergence and consensus of the DRGFMD algorithm.}\label{config1}
\end{figure}

\begin{figure}[htbp]
\centering
\includegraphics[width=2in]{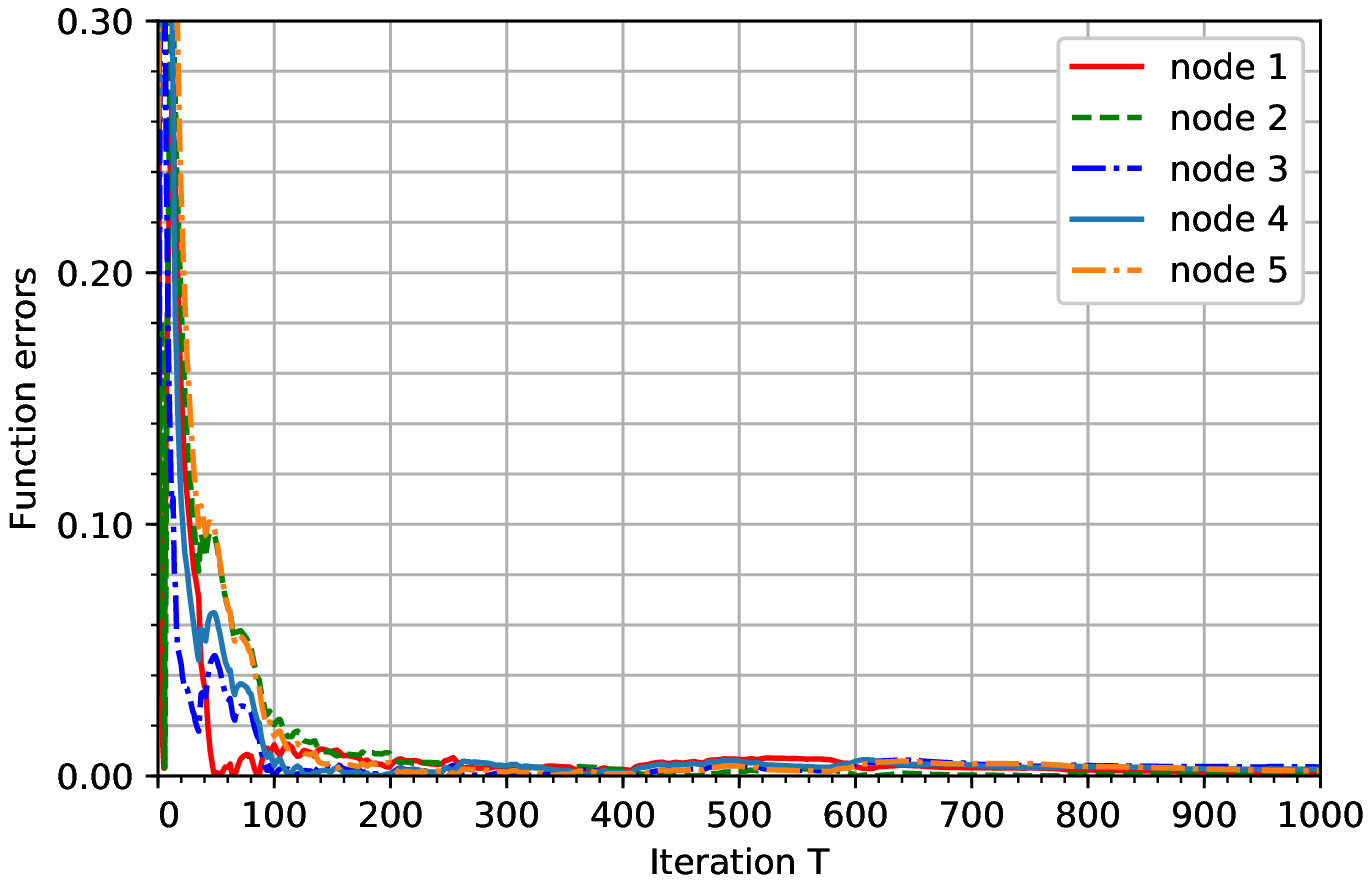}
\caption{Convergence and consensus of the DRGFMD-WA algorithm.}\label{config2}
\end{figure}

\begin{figure}[htbp]
\centering
\includegraphics[width=2in]{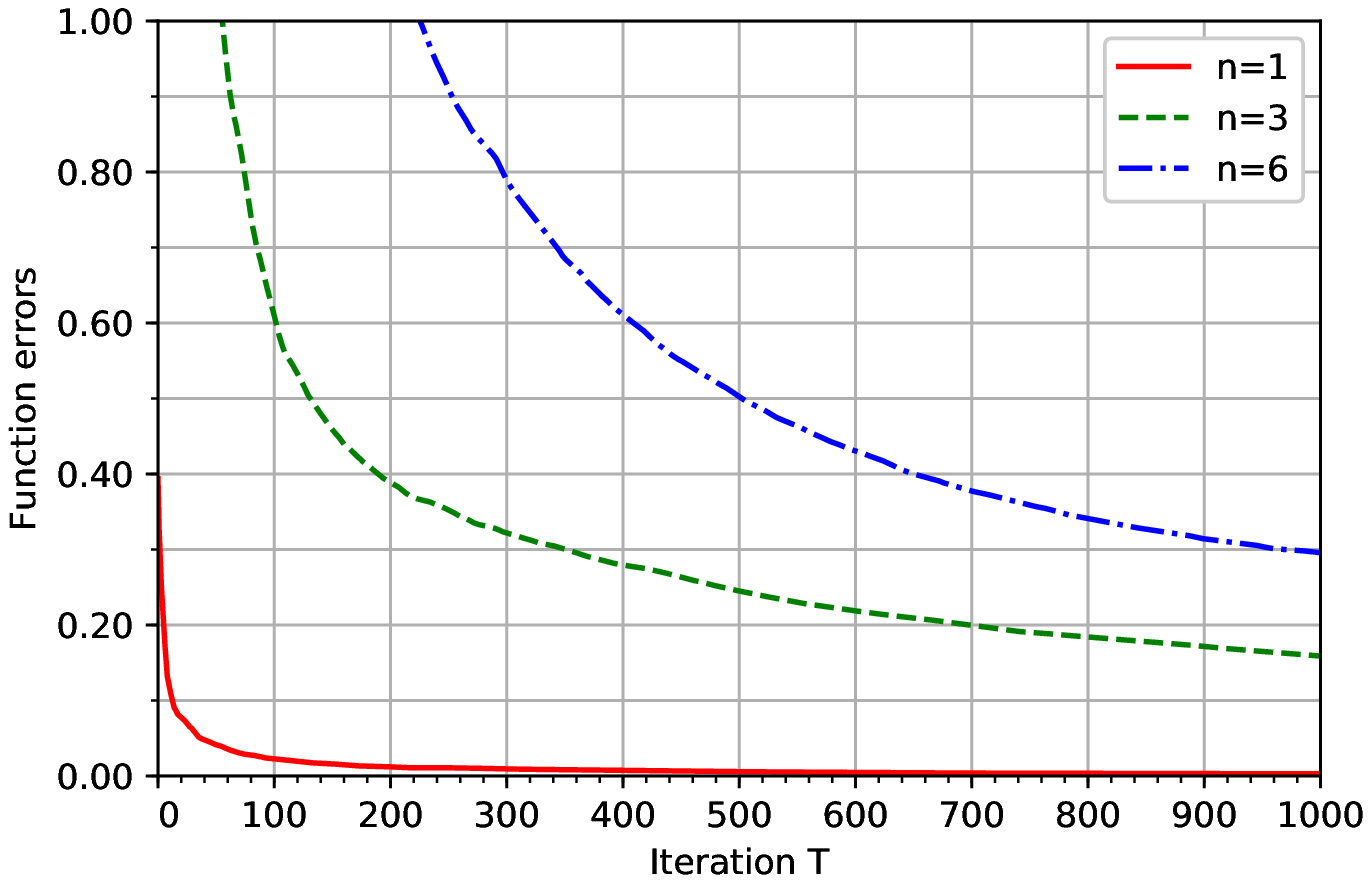}
\caption{Influence of dimension of agents' estimates on the convergence of DRGFMD algorithm. }\label{config3}
\end{figure}

\begin{figure}[htbp]
\centering
\includegraphics[width=2in]{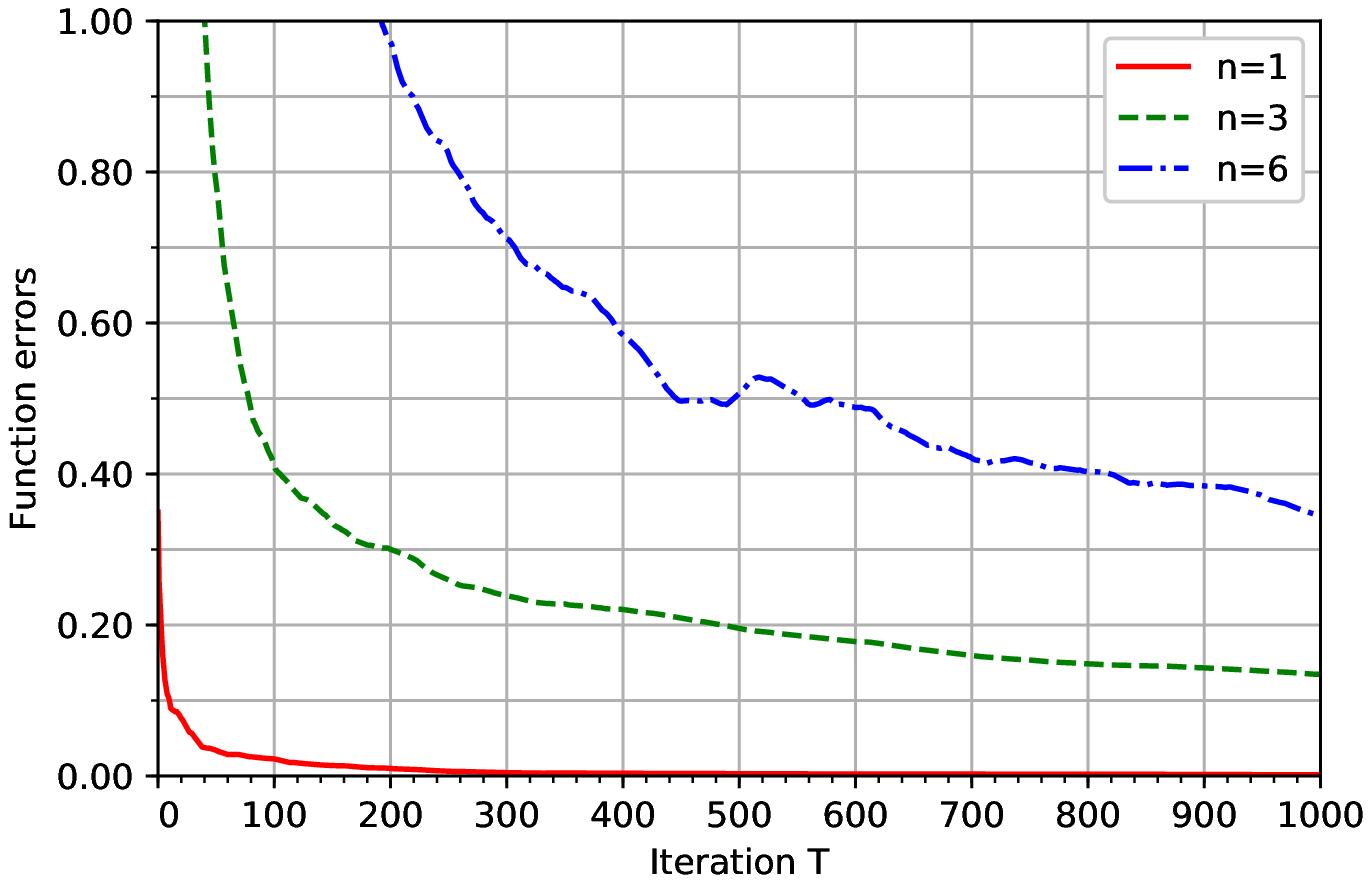}
\caption{Influence of dimension of agents' estimates on the convergence of DRGFMD-WA algorithm.}\label{config4}
\end{figure}

\begin{figure}[htbp]
\centering
\includegraphics[width=2in]{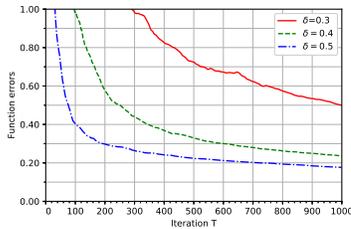}
\caption{Influence of stepsize on the convergence of DRGFMD-WA algorithm.}\label{config5}
\end{figure}

In this section, several descriptions of the DRGFMD algorithm and DRGFMD-WA algorithm are given by providing a simulation example. Specifically, the DRGFMD algorithm and DRGFMD-WA algorithm are utilized to analyze the Nesterov nonsmooth test problem given by
\begin{eqnarray}
\nono\min_{x\in X} \sum_{i=1}^N c_i(|[x]_1-1|+\sum_{k=1}^{n-1}|1+[x]_{k+1}-2[x]_{k}|)
\end{eqnarray}
in constraint set $X=\{x\in \mathbb R^n: \sum_{d=1}^n[x]_d=1,[x]_d\geq0\}$, and $c_i\in \mathbb R^{+}$ is the data known only to node $i$. The random graph with nodes $N=5$ which is generated in a manner of Xiao ad Boyd (\cite{s6}) is considered. The gradient-free random sequence $\xi_i^t$ is generated in an independent and identically distributed way from Gaussian normal distribution $N(0,0.5I_{n\times n})$ for all $i$. $\phi(x)=\sum_{d=1}^n[x]_d\ln[x]_d$ is chosen as the distance-generating function of the proposed algorithms. In following simulations, the DRGFMD algorithm and the DRGFMD-WA algorithm are used as the trial objects. The first two simulation results use $n=1$ and $\mu_i=10^{-4}$ to give an obvious description on convergence and consensus behavior among agents as tests in Figure \ref{config1} and Figure \ref{config2}. The results coincide with the fact that two algorithms achieve the same convergence rate under same convexity condition on objective functions and constraint conditions on $X$. In what follows, the  influence of the dimension of the decision space on the convergence of DRGFMD algorithm and DRGFMD-WA algorithm is considered.  Simulations of $n=1$, $n=3$ and $n=6$ are investigated to reflect the distinct difference of the convergence speed with different dimension. In each of these cases, initial data $x_i^0=(1/n, 1/n,...,1/n)^T$ and $\mu_i=10^{-4}$ are used. The simulation results on dimension influence are based on the average of 30 independent trials. Both of Figures \ref{config3} and \ref{config4} reveal that the convergence is faster with smaller dimension. That is to say, when $n$ becomes larger, it will take more iterations for DRGFMD algorithm and DRGFMD-WA algorithm to reach the same accuracy as the one with smaller $n$. The simulations on the influence of dimension  on the convergence is not accidental, since they are in compliance with the convergence results of the paper. After simulating the dimension influence, the influence of stepsize on the DRGFMD-WA algorithm is considered in Figure \ref{config5}. The simulation uses $n=3$, $\mu_i=10^{-4}$ and $x_i^0=(1/3, 1/3, 1/3 )^T$, stepsize $\aaa_t=1/(t+1)^{\delta}$ ($\delta\in (0,1/2]$). The simulation result of each stepsize is an average of 30 independent trials. In the trial, three different values $\delta=0.3$, $\delta=0.4$, $\delta=0.5$ are investigated. The selected $\delta$s are linear 1/10-increasing values. Figure \ref{config5} shows that with $\delta$ getting larger, the convergence becomes faster which is as expected. However, there is an obviously bigger gap between the case $\delta=0.3$ and $\delta=0.4$, which indicates that, with $\delta$ getting smaller, the degree that the convergence of DRGFMD-WA gets slower becomes obviously greater. The next simulation gives a comparison among DRGFMD, DRGFMD-WA and previous distributed gradient-free projection method (DGFP) with approximating sequence $\widehat{x_i}^T=\sum_{t=1}^T\aaa_t x_i^t/\sum_{t=1}^T\aaa_t$ in \cite{gf}. This simulation uses $n=2$, $\mu_i=10^{-4}$ and initial data $x_i^0=(1/2,1/2)^T$ in each algorithm, the comparison result is based on the average of 30 independent trials. The result shows DRGFMD and DRGFMD-WA (accuracy around $10^{-2}$) are much more efficient than the previous DGFP algorithm (accuracy around $10^{-1}$).

\section{Conclusions and discussions}
In this paper, both convex and strongly convex constrained distributed optimization problem are considered by developing a distributed randomized gradient-free mirror descent method. To implement the DRGFMD method, the gradient and subgradient information of objective functions is not necessary to be known. The convergence rates of the DRGFMD algorithm are considered under two types of conditions. A decentralized reciprocal weighted average approximating sequence is first investigated in DRGFMD framework and a class of corresponding convergence rates are achieved. Finally the simulation results are presented to illustrate the convergence behavior in several aspects.

The work in this paper opens a few future research directions. One is the construction of the decentralized weighted average approximating sequence. The idea and technique of the construction can be generally used in large amount of distributed algorithms. Further,  using the proposed scheme of decentralized weighted averaging, the convergence rates of distributed algorithms can be investigated, and potential improvement is expected.  Also, further research on the algorithm can be explored (i) by utilizing appropriate diminishing smoothing parameter to eliminate the effect on error bound; (ii) by constructing appropriate distributed zeroth-order oracles to reduce the large dimension influence on convergence. Moreover, other possible  application directions are to extend the  proposed algorithm to distributed online optimization, and to apply the proposed zeroth-order method to distributed nonconvex optimization problem.


\vspace{-6 mm}
\vspace{-6 mm}
\vspace{-5 mm}
\vspace{-4 mm}
\vspace{-4 mm}
\begin{figure}[htbp]
\centering
\includegraphics[width=2in]{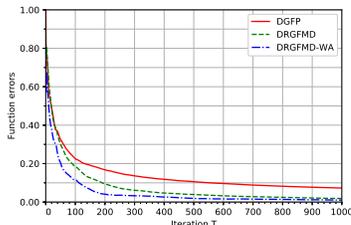}
\caption{Comparison among DRGFMD algorithm, DRGFMD-WA algorithm and DGFP algorithm.}\label{config7}
\end{figure}

\end{CJK*}
\end{document}